\documentclass[3p,11pt]{elsarticle}

\usepackage{amssymb}
\usepackage{amsmath}
\usepackage{amsfonts}
\usepackage{indentfirst}
\usepackage{graphicx}
\usepackage[english,francais]{babel}

\usepackage{tikz}
\usepackage{pgfplots}

\newtheorem{rem}{Remark}

\newcommand{\vps}{\varepsilon}
\newcommand   {\ep}  {\varepsilon}

\newcommand {\noame} {\noalign{\medskip}}
\newcommand {\dis} {\displaystyle}

\newcommand{\beq}{\begin{equation}}
\newcommand{\eeq}{\end{equation}}

\newcommand{\LV}{L^2_{x_1}(0,1;V_Z)}
\newcommand{\LL}{L^2_{x_1}(0,1;L^2_Z(0,1))}

\newcommand{\RR}{\mathbf{R}}

\newcommand{\wra}{\rightharpoonup}
\newcommand{\embed}{\hookrightarrow}

\newcommand{\ub}{\mathbf{u}}

\newcommand{\wb}{\mathbf{w}}

\newcommand{\ueps}{\ub_{\vps}}
\newcommand{\weps}{\wb_{\vps}}

\renewcommand{\d}{\partial}
\newcommand{\uzero}{u_1^0}

\newcommand{\un}{u_1^n}
\newcommand{\wn}{w_2^n}
\newcommand{\unp}{u_1^{n+1}}

\newcommand{\wnp}{w_2^{n+1}}

\newcommand{\unptilde}{\tilde{u}_1^{n+1}}
\newcommand{\pnp}{p^{n+1}}

\newcommand{\psibar}{\bar{\psi}}

\newcommand{\be}{\begin{equation}}
\newcommand{\ee}{\end{equation}}
\newcommand{\bestar}{\begin{equation*}}
\newcommand{\eestar}{\end{equation*}}

\begin{document}

\begin{frontmatter}


\selectlanguage{english}
\title{Effects of rough boundary and nonzero boundary conditions on the lubrication process with micropolar fluid}

\selectlanguage{english}
\author[authorlabel1]{Matthieu Bonnivard}
\ead{bonnivard@ljll.univ-paris-diderot.fr}
\author[authorlabel2]{Igor Pa\v{z}anin}
\ead{pazanin@math.hr}
\author[authorlabel3]{Francisco Javier Su\'arez-Grau}
\ead{fjsgrau@us.es}
\address[authorlabel1]{Univ. Paris Diderot, Sorbonne Paris Cit\'e, Laboratoire Jacques-Louis Lions, UMR 7598, UPMC, CNRS, F-75205 Paris, France}
\address[authorlabel2]{Department of Mathematics, Faculty of Science, University of Zagreb,
Bijeni\v{c}ka 30, 10000 Zagreb, Croatia}
\address[authorlabel3]{Departamento de Ecuaciones Diferenciales y An\'alisis Num\'erico, Facultad de Matem\'aticas, Universidad de Sevilla, Avenida Reina Mercedes S/N, 41012 Sevilla, Spain}

\medskip

\begin{abstract}
The lubrication theory is mostly concerned with the behavior of a lubricant flowing through a narrow gap. Motivated by the experimental findings from the tribology literature, we take the lubricant to be micropolar fluid and study its behavior in a thin domain with rough boundary. Instead of considering (commonly used) simple zero boundary condition, we impose physically relevant (nonzero) boundary condition for microrotation and perform asymptotic analysis of the corresponding 3D boundary value problem. We formally derive a simplified mathematical model acknowledging the roughness-induced effects and the effects of the nonzero boundary conditions on the macroscopic flow. Using the obtained asymptotic model, we study numerically the influence of the specific rugosity profile on the performance of a linear slider bearing. The numerical results clearly indicate that the use of the rough surfaces may contribute to enhance the mechanical performance of such device.

\vskip 0.5\baselineskip

\keyword{lubrication process, micropolar fluid, nonzero boundary conditions, rough boundary, asymptotic analysis, numerical results.}
}
\end{abstract}
\end{frontmatter}

\selectlanguage{english}

\section{Introduction}
In this paper we investigate the flow of a viscous fluid (acting as a lubricant) between two solid surfaces in  relative motion. The lower surface is assumed to be plane, while the upper is described by a given shape function. Such situation appears naturally in numerous industrial and engineering applications, in particular those consisting of moving machine parts (see e.g.~\cite{Szeri}). The mathematical models for describing the motion of the lubricant usually result from the simplification of the geometry of the lubricant film, i.e.~its thickness. Using the film thickness as a small parameter, a simple asymptotic approximation can be easily derived providing a well-known Reynolds equation for the pressure of the fluid. Formal derivation goes back to the 19th century and the celebrated work of Reynolds \cite{Reynolds}. The justification of this approximation, namely the proof that it can be obtained as the limit of the Stokes system (as thickness tends to zero) is provided in \cite{Bayada1} for a Newtonian flow between two plain surfaces. More recent results on the lubrication with a Newtonian fluid can be found in \cite{Archiv,Wil,Ja1}.

The experimental observations from the tribology literature (see e.g.~\cite{exp1,exp2,exp3}) indicate that the fluid's internal structure should not be ignored in the modelling, especially when the gap between the moving surfaces is very small. Among various non-Newtonian models, the model of micropolar fluid (proposed by Eringen \cite{Eringen} in 60's) turns out to be the most appropriate since it acknowledges the effects of the local structure and micro-motions of the fluid elements. Physically, micropolar fluids consist in a large number of small spherical particles uniformly dispersed in a viscous medium. Assuming that the particles are rigid and ignoring their deformations, the related mathematical model expresses the balance of momentum, mass and angular momentum. A new unknown function called microrotation (i.e.~the angular velocity field of rotation of particles) is added to the usual velocity and pressure fields. Consequently, Navier-Stokes equations become coupled with a new vector equation coming from the conservation of angular momentum with four microrotation viscosities introduced. Being able to describe numerous real fluids better than the classical (Newtonian) model, micropolar fluid model has been extensively studied in recent years (see e.g.~\cite{Panasenko1,Panasenko2,Ja2,Ja3,Ja4}.)

Engineering practice also stresses the interest of studying the effects of very small domain irregularities on a thin film flow. That means that the lower surface in our setting is assumed to be perfectly smooth, while the upper is rough with roughness described by a given (periodic) function. Assuming specific rugosity profile in which the size of the rugosities is smaller than the film thickness, Bresch and co-authors \cite{Bresch} derived the explicit correction of the Reynolds approximation by the roughness-induced term. Though specific, such roughness pattern is physically relevant and realistic (see e.g.~\cite{fiz}), and results in modifying the Reynolds equation at the main order\footnote{The usual assumption that the size of the rugosities is (at least) of the same order as the film thickness leads to the effective models in the form of the classical Reynolds equation. In order to detect the roughness-induced effects on the macroscopic flow, one has to compute the lower-order terms (see e.g.~\cite{Bayada2,Benh}).}.

In the above mentioned papers (see also \cite{Chupin,Ja5}), the lubricant is assumed to be a classical Newtonian fluid. So far, micropolar fluid film lubrication has been addressed mostly in a simple thin domain with no roughness introduced (see e.g.~\cite{Sinha2,Bayada3,Ja6}). However, recently, the roughness effects on a thin film flow have been studied as well and new mathematical models have been proposed. Namely, in \cite{Boukrouche}, the authors consider micropolar flow in a 2D domain assuming the roughness is of the same small order as the film thickness. Taking zero boundary condition for microrotation, the effective system is derived using two-scale convergence method. In \cite{Ja7}, a generalized version of the Reynolds equation is derived in the case of a 3D domain with roughness pattern proposed in \cite{Bresch}. Again, simple zero boundary condition is imposed on the microrotation, implying that the fluid elements cannot rotate on the fluid-solid interface. Though the physical interpretation of such boundary condition is clearly doubtful, it has been commonly used throughout the literature. The aim of this paper is to consider other (physically relevant) type of boundary condition for microrotation originally proposed in \cite{Bessonov1,Bessonov2}. This new type of boundary condition is based on the concept of boundary viscosity\footnote{Experimental observations suggest that, in a thin film flow of a non-Newtonian fluid, the chemical interactions between solid surface and the lubricant should not be ignored. It turns out that this phenomenon can be taken into account by introducing the so called boundary viscosity, which differs from Newtonian and microrotation viscosities.} linking the value of the microrotation with the rotation of the velocity on the fluid-solid interface. In such setting, standard variation principle suggests that the no-slip boundary condition for the velocity needs to be replaced as well, namely by the condition allowing slippage of the wall. This particular feature makes the asymptotic analysis of the corresponding boundary value problem more demanding than in the standard case with zero boundary conditions.

The 2-dimensional micropolar flow associated to the above described (nonstandard) boundary conditions in a thin domain without roughness has been investigated in \cite{Bayada4} (see also \cite{Ja8} for the 3D flow). A generalized version of the Reynolds equation has been rigorously derived in a critical case where one of the non-Newtonian characteristic parameters has specific (small) order of magnitude. Here we intend to study a similar problem but in a rough domain, with roughness described by a periodic function with period of order $\ep^{2}$ (the domain's thickness is assumed to be of order $\ep$). Employing a formal asymptotic analysis with respect to $\ep$, in Sec.~3 we derive a strongly coupled effective system, which clearly indicates the influence of the rough boundary and nonzero boundary conditions on the macroscopic flow. Moreover, in Sec.~4 we conduct numerical simulations based on the  asymptotic model obtained formally, for a particular lubrication device - a linear slider bearing with rough inclined surface. We develop a numerical strategy successfully dealing with the technical difficulties arising due to the effect of the roughness and the nonstandard boundary conditions. As a result, the simulations show that the considered rugosity profile may, in fact, contribute to enhance the mechanical performance of the lubrication device. In view of that, we believe that the result presented in this paper could be instrumental for improving the known engineering practice.

\section{Position of the problem}
The domain under consideration has the following form:
$$\widehat \Omega_\ep=\left\{(\widehat x,\widehat z)\in\mathbf{R}^2\times \RR\ :\  \widehat x\in \omega,\quad 0<\widehat z<\widehat{h}_{\varepsilon}(\widehat{x})\right\}\,.$$
Here $\omega\subset \mathbf{R}^2$ is an open subset with smooth boundary, and $\hat{h}_{\vps}(\widehat{x})$ denotes the height of the domain given by
\begin{eqnarray}\label{rug-hat}
\widehat h_\ep(\widehat x) = c\,h_1\left(\frac{\widehat x}{L}\right) + \frac{c^2}{L}\,h_2\left(\frac{L}{c^2}\widehat x \right)\,.
\end{eqnarray}
In this above expression, $h_2$ is a regular $\RR$-valued function, defined on $\RR^2$,  periodic with period $[0,1]^2$ and with zero average, that models the roughness of the upper surface.
Moreover, L denotes the characteristic length of $\omega$, whereas $c$ stands for the maximum distance between the surfaces. We consider the ratio $\ep=\frac{c}{L}$ as the small parameter.

Let $\widehat{\Gamma}^1_\ep$, $\widehat{\Gamma}_\ep^0$ and $\widehat{\Gamma}_{\ep}^\ell$ denote the upper, lower and lateral boundaries of $\widehat{\Omega}_\varepsilon$, namely:
$$\begin{array}{l}
\dis \widehat{\Gamma}^1_\ep=\left\{(\widehat x,\widehat z)\in\RR^2\times \RR\ :\  \widehat x\in \omega,\quad \widehat z=\widehat{h}_{\varepsilon}(x)\right\},\\
\noame\dis
\widehat{\Gamma}_\ep^0=\left\{(\widehat x,\widehat z)\in\RR^2\times \RR\ :\  \widehat x\in \omega,\quad \widehat z=0\right\},\\
\noame\dis
\widehat{\Gamma}_{\ep}^\ell=\partial\widehat{\Omega}_\varepsilon-(\widehat{\Gamma}_\ep^1\cup \widehat{\Gamma}_{\ep}^0).
\end{array}$$
The micropolar fluid flow is described by the following equations\footnote{Note that here we consider the linearized equations, as in \cite{Bayada4}. Taking into account the application we want to model (lubrication with micropolar fluid), it is reasonable to assume a small Reynolds number and omit the inertial terms in momentum equations.} expressing the balance of momentum, mass and angular momentum (see \cite{Lukas}):
\begin{eqnarray}
 &&-(\nu+\nu_r)\,\Delta\widehat{\mathbf{u}}_\ep + \nabla \widehat{p}_\ep =2\nu_r\mbox{rot}\,\widehat{\mathbf{w}}_{\ep}\,\,\,\,\,\,\,\,\,\mbox{in}\,\,\,\,\,\,\widehat{\Omega}_\varepsilon\;,\label{1,1-hat}\\[0.2cm]
 &&\mbox{div}\,\widehat{\mathbf{u}}_\ep =0\,\,\,\,\,\,\,\,\,\mbox{in}\,\,\,\,\,\,\widehat{\Omega}_\varepsilon\;,\label{21,1-hat}\\[0.2cm]
 &&-(c_a+c_d)\,\Delta\,\widehat{\mathbf{w}}_\ep+4\nu_r\widehat{\mathbf{w}}_{\ep}=2\nu_r\mbox{rot}\,\widehat{\mathbf{u}}_{\ep}\,\,\,\,\,\,\,\,\,\mbox{in}\,\,\,\,\,\,\widehat{\Omega}_\varepsilon\;.\label{31,1-hat}
\end{eqnarray}
In the above system, velocity $\widehat{{\bf u}}_\ep$, pressure $\widehat{p}_\ep$ and microrotation $\widehat{{\bf w}}_\ep$ are unknown. $\nu$ is the Newtonian viscosity, while $\nu_r$, $c_a$, $c_d$ are microrotation viscosities appearing due to the asymmetry of the stress tensor. All viscosity coefficients are positive constants.

As discussed in the Introduction, we impose the following boundary conditions:
\begin{eqnarray}
&&\widehat{\mathbf{u}}_\varepsilon = 0,\quad \widehat{\mathbf{w}}_\varepsilon = 0 \quad\hbox{    on    }\quad\widehat{\Gamma}_\ep^1,\label{bc1}\\[0.2cm]
&&\widehat{\mathbf{u}}_\varepsilon = \widehat{\mathbf{g}},\quad \widehat{\mathbf{w}}_\varepsilon = 0\quad \hbox{    on    }\quad\widehat{\Gamma}_\ep^\ell,\label{bc2}\\[0.2cm]
&&\widehat{\mathbf{u}}_\varepsilon\cdot \mathbf{k} = 0,\quad \widehat{\mathbf{w}}_\varepsilon\cdot \mathbf{k} = 0\quad \hbox{    on     }\quad\widehat{\Gamma}_\ep^0,\label{bc3}\\[0.2cm]
&&\widehat{\mathbf{w}}_\varepsilon\times \mathbf{k} ={\alpha\over 2}({\rm rot}\widehat{\mathbf{u}}_\ep)\times \mathbf{k}\quad \hbox{    on     }\quad\widehat{\Gamma}_\ep^0,\label{bc4}\\[0.2cm]
&&({\rm rot}\widehat{\mathbf{w}}_\ep)\times \mathbf{k}={2\nu_r\over c_a+c_d}\beta(\widehat{\mathbf{u}}_\ep-\widehat{\mathbf{s}})\times \mathbf{k}
\quad \hbox{    on     }\quad\widehat{\Gamma}_\ep^0.\label{bc5}
\end{eqnarray}
Here and in the sequel $( \textbf{i},\textbf{j},\textbf{k})$ denotes the standard
Cartesian basis. Notice that the usual (Dirichlet) boundary conditions (\ref{bc1})-(\ref{bc2}) for the velocity and microrotation are prescribed on $\Gamma_\ep^1\cup \Gamma_\ep^\ell$. However, on the lower part $\Gamma_\ep^0$ (corresponding to the moving boundary), new type of boundary conditions (\ref{bc4})-(\ref{bc5}) are imposed\footnote{The impermeability of the wall leads to (\ref{bc3}).}. The coefficient $\alpha>0$ appearing in (\ref{bc4}) characterizes the microrotation on the solid surfaces and it is given by $\alpha={\nu+\nu_r-\nu_b\over \nu_r}$ (see \cite{Bessonov1,Bessonov2}). Here $\nu_b$ denotes a new viscosity coefficient introduced in the vicinity of the surface, and different from $\nu$ and $\nu_r$. Finally, in order to obtain a well-posed variational formulation, it is straightforward to confirm (see e.g.~\cite{Bayada4}) that boundary condition (\ref{bc5}) needs to be introduced. The coefficient $\beta$ in (\ref{bc5}) allows the control of the slippage at the wall when the value $\widehat{\mathbf{u}}_\ep-\widehat{\mathbf{s}}$ is not zero ($\widehat{\mathbf{s}}$ is a given velocity of the wall).

It has been observed (see e.g.~\cite{Bayada3,Bayada4,Ja9}) that the magnitude of the viscosity coefficients appearing in (\ref{1,1-hat})-(\ref{31,1-hat}) may influence the effective flow. Thus, it is reasonable to work with the system written in a non-dimensional form. In view of that, we introduce:
\begin{eqnarray}
&&x={\widehat x\over L},\quad z={\widehat z\over L},\quad h_\varepsilon={\widehat h_\varepsilon\over L},\nonumber\\[0.2cm]
&&\mathbf{u}_\ep={\widehat{\mathbf{u}}_\varepsilon\over s_0},\quad p_\varepsilon={L\over s_0(\nu+\nu_r)}\widehat p_\varepsilon,\quad \mathbf{w}_\ep={L\over s_0}\widehat{ \mathbf{w}}_\varepsilon,\quad \mathbf{g}={\widehat{\mathbf{g}}\over s_0},\quad \mathbf{s}={\widehat{\mathbf{s}}\over s_0},\nonumber\\[0.2cm]
&&N^2={\nu_r\over \nu+\nu_r},\quad R_M={c_a+c_d\over \nu+\nu_r}{1\over L^2}\,,\label{nonparam}
\end{eqnarray}
where $s_0$ stands for the characteristic velocity of the moving surface. Dimensionless (non-Newtonian) parameter $N^{2}$ characterizes the coupling between the equations for the velocity and microrotation and it is of order $\mathcal{O}(1)$. The second dimensionless parameter, denoted by $R_M$
is, in fact, related to the characteristic length of the microrotation effects and will be compared with small parameter $\ep$.

In view of the above change, the fluid domain becomes
$$\Omega_\ep=\left\{(x,z)\,:\, x\in \omega/L,\quad 0<z<h_\ep(x)\right\}\,.$$
The height of the domain is now given by (see (\ref{rug-hat}))
\begin{eqnarray}
h_\varepsilon(x)=\varepsilon\,  h_1(x)+\varepsilon^2\, h_2\left({x\over \varepsilon^2}\right)\,.\label{rugprofile}\end{eqnarray}
The functions $h_1,h_2$ are assumed to be regular: the positive function $h_1\in H^2(\omega/L)$ represents the main order part, while the $\mathbb{T}^2$-periodic function $h_2\in H^2(\mathbb{T}^2) $ (with $0$ as average in $\mathbb{T}^2$, $\mathbb{T}^2$  is the torus of dimension $2$) describes the oscillating part of the roughness (see Fig.~1). This rugosity profile (\ref{rugprofile}) has been addressed in \cite{Bresch} for classical, Newtonian flow and in \cite{Ja7} for micropolar flow associated to zero boundary conditions.
\begin{figure}[h!]
\centerline{
\includegraphics[width=12cm]{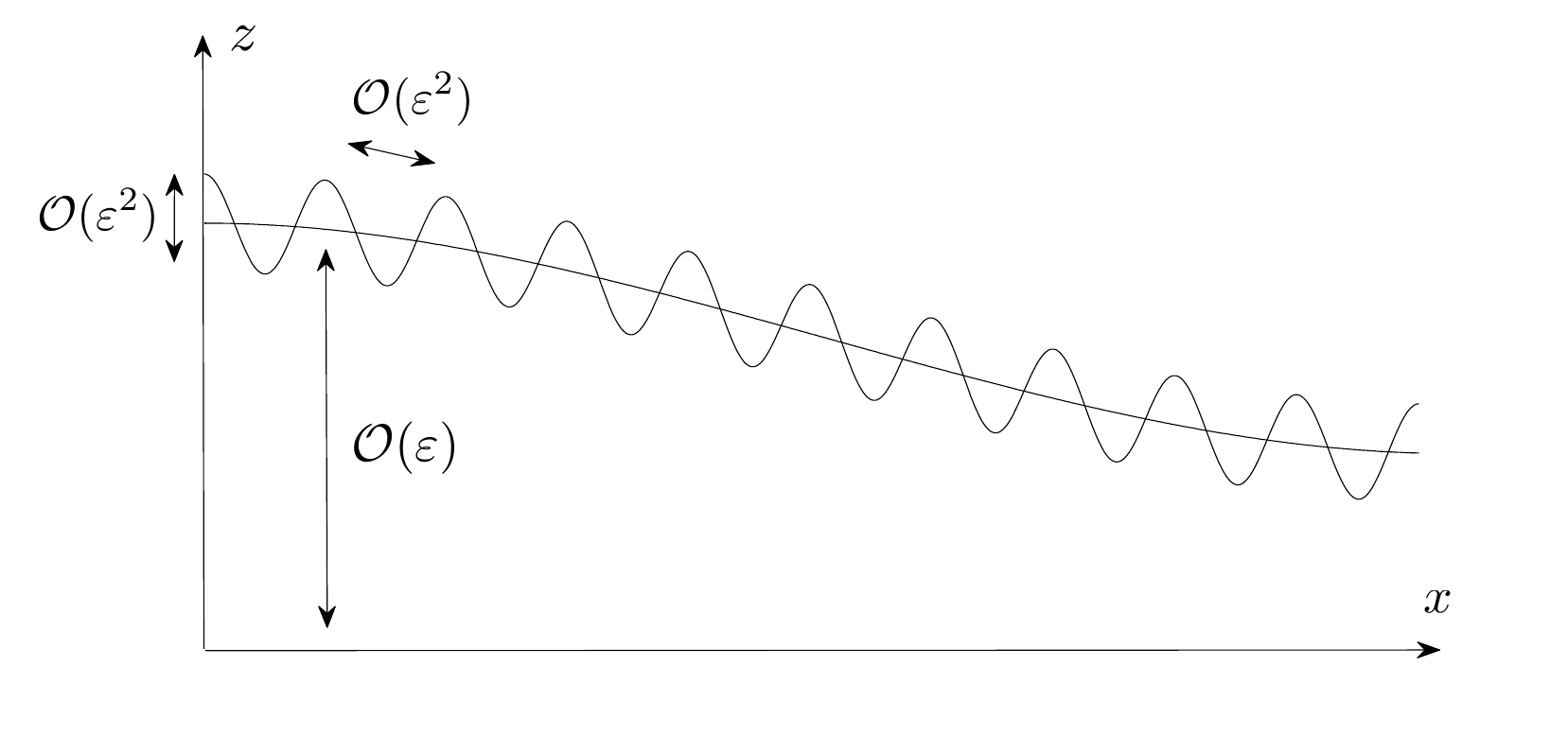}}
\caption{The different scales related to the domain.}
\end{figure}

The flow equations (\ref{1,1-hat})-(\ref{31,1-hat}) now have the following form
\begin{eqnarray}
 &&-\,\Delta\mathbf{u}_\ep + \nabla p_\ep =2 N^2\mbox{rot}\,\mathbf{w}_{\ep}\,\,\,\,\,\,\,\,\,\mbox{in}\,\,\,\,\,\,\Omega_\varepsilon\;,\label{1,1}\\[0.2cm]
 &&\mbox{div}\,\mathbf{u}_\ep =0\;\label{21,1}\,\,\,\,\,\,\,\,\,\mbox{in}\,\,\,\,\,\,\Omega_\varepsilon\,,\\[0.2cm]
 &&-R_M\,\Delta\,\mathbf{w}_\ep+4N^2\mathbf{w}_{\ep}=2N^2\mbox{rot}\,\mathbf{u}_{\ep}\,\,\,\,\,\,\,\,\,\mbox{in}\,\,\,\,\,\,\Omega_\varepsilon\;.\label{31,1}
\end{eqnarray}
Let us homogenize the boundary condition (\ref{bc2})$_1$. Assuming the existence of ${\bf g}\in H^{1\over 2}(\Gamma_{\ep}^\ell)$ such that $\int_{\Gamma_\ep^\ell}\bf g\cdot \mathbf{n}\, d\sigma =0$, one can easily construct (see e.g.~\cite{Bayada4}) a lift function ${\bf J}_\ep\in H^1(\Omega_\ep)$ such that
\begin{equation}\label{function-J}
{\rm div}\mathbf{J}_\ep=0\ \hbox{ in }\Omega_\ep,\quad \mathbf{J}_\ep=0\ \hbox{ on }\Gamma_\ep^1, \quad \mathbf{J}_\ep=\mathbf{g}\ \hbox{ on }\Gamma_\ep^\ell,\quad J_3^\ep=0\ \hbox{ on }\Gamma_\ep^0.
\end{equation}
Now we introduce
\begin{equation}\label{v-u-J}{\bf v}_\ep={\bf u}_\ep-{\bf J}_\ep\,.
\end{equation}
In view of that, the system (\ref{1,1})-(\ref{31,1}) can be rewritten as
\begin{eqnarray}
 &&-\Delta(\mathbf{v}_\ep+\mathbf{J}_\ep) + \nabla p_\ep =2N^2\mbox{rot}\,\mathbf{w}_{\ep}\,\,\,\,\,\,\,\,\,\mbox{in}\,\,\,\,\,\,\Omega_\varepsilon\;,\label{1,1-J}\\[0.2cm]
 &&\mbox{div}\,\mathbf{v}_\ep =0\,\,\,\,\,\,\,\,\,\mbox{in}\,\,\,\,\,\,\Omega_\varepsilon\;,\label{21,1-J}\\[0.2cm]
  &&-R_M\Delta\,\mathbf{w}_\ep+4N^2\mathbf{w}_{\ep}=2N^2(\mbox{rot}\,\mathbf{v}_{\ep} +\mbox{rot}\,\mathbf{J}_{\ep} )\,\,\,\,\,\,\,\,\,\mbox{in}\,\,\,\,\,\,\Omega_\varepsilon\;.\label{31,1-J}
\end{eqnarray}
The corresponding boundary conditions read
\begin{eqnarray}
&&\mathbf{v}_\varepsilon = 0,\quad \mathbf{w}_\varepsilon = 0 \quad\hbox{    on    }\quad\Gamma_\ep^1\cup\Gamma_\ep^\ell,\label{bc-1}\\[0.2cm]
&&\mathbf{v}_\varepsilon\cdot \mathbf{k} = 0,\quad \mathbf{w}_\varepsilon\cdot \mathbf{k} = 0\quad \hbox{    on     }\quad\Gamma_\ep^0,\label{bc-3}\\[0.2cm]
&&\mathbf{w}_\varepsilon\times \mathbf{k} ={\alpha\over 2}(({\rm rot}(\mathbf{v}_\ep+\mathbf{J}_\ep))\times \mathbf{k}\quad \hbox{    on     }\quad\Gamma_\ep^0,\label{bc-4}\\[0.2cm]
&&R_M({\rm rot}\mathbf{w}_\ep)\times \mathbf{k}=2N^2\beta(\mathbf{v}_\ep+\mathbf{f}-\mathbf{s})\times \mathbf{k}\quad \hbox{    on     }\quad\Gamma_\ep^0.\label{bc-5}
\end{eqnarray}

\begin{rem}
It should be observed that only the normal component of the velocity is known on $\Gamma_0^\ep$, while the tangential component is not given (see (\ref{bc-3})$_1$). Nevertheless, we can choose an artificial value $\mathbf{f}=(f_1,f_2)$ of the velocity on $\Gamma_0^\ep$ appearing in $(\ref{bc-5})$. We choose it in a way such that function $\mathbf{d}\in H^{1\over 2}(\partial\Omega^\ep)$ defined by
$$\mathbf{d}=\left\{\begin{array}{l}
0\quad\hbox{ on }\quad \Gamma_\ep^1,\\
\noame
\mathbf{g}\quad\hbox{ on }\quad \Gamma_\ep^\ell,\\
\noame
(\mathbf{f},0)\quad\hbox{ on }\quad \Gamma_\ep^0
\end{array}\right.$$
satisfies $\int_{\partial\Omega_\ep}\mathbf{d}\cdot\mathbf{n}=0$. Consequently, we are in position to construct the lift function $\mathbf{J}_\ep$ as in (\ref{function-J}).
\end{rem}

The well-posedness of the above problem can be established using the same arguments as in \cite{Bayada4}. In the present paper the aim is to derive the macroscopic law describing the effective flow in $\Omega_\ep$ by using asymptotic analysis with respect to the small parameter $\ep$. In particular, we shall focus on detecting the roughness-induced effects and the effects of nonzero boundary conditions.


\section{Asymptotic analysis}
As indicated before, different asymptotic behaviors of the flow may be deduced depending on the order of magnitude of the viscosity coefficients. Indeed, if we compare the characteristic number $R_M$ defined by (\ref{nonparam})$_2$ and appearing in the equation (\ref{31,1-J}) with small parameter $\ep$, three different asymptotic situations can be formally identified (see e.g.~\cite{Ja9}). The most interesting one is, of course, the one leading to a strong coupling at a main order, namely when
\begin{eqnarray}
&&R_M=\ep^2 R_c\,,\,\,\,\,\,\,\,\,R_c=\mathcal{O}(1)\,.
\end{eqnarray}
Thus, we will proceed with our analysis assuming the above scaling of $R_M$ with respect to $\ep$.

\subsection{Rescaling}
The first step in the asymptotic analysis is to rewrite the governing problem in $\ep$-independent domain. This should be done in a way so the rescaling procedure acknowledges the oscillating behavior of our thin domain. Following \cite{Bresch}, we introduce a fast variable
$$X=\frac{x}{\ep^2}\,.$$
In view of that, the height $h_\ep$ becomes
$$h(x,X)=\ep h_1(x)+\ep^2 h_2(X)\,.$$
We also introduce new (vertical) variable $Z=\frac{z}{h(x,X)}$ and, correspondingly, the new unknown functions:
$$\mathbf{u}_\ep(x,z)=\mathbf{\tilde{u}}(x,X,Z),\,\,\,\, \mathbf{w}_\ep(x,z)=\mathbf{\tilde{w}}(x,X,Z),\,\,\,p_\ep(x,z)=p(x,X,Z)\,.$$
In the following, we set
$$\mathbf{\tilde{u}}=(\mathbf{v},v_3)\in\RR^{2}\times\RR,\quad \mathbf{\tilde{w}}=(\mathbf{w},w_3)\in\RR^{2}\times\RR\,.$$
The lift $\mathbf{J}_\ep=(J_{1,\ep},J_{2,\ep},J_{3,\ep})$ (introduced in the previous section) can be constructed from $\overline{\mathbf{J}}=(J_1,J_2,J_3)$ defined in $\ep$-independent domain $\Omega$ such that ${\rm div}\overline{\mathbf{J}}=0$ by putting:
$$J_{i,\ep}(x,z)=J_i\left(x,{z\over h_\ep(x)}\right),\ i=1,2,\quad J_{3,\ep}(x,z)=h_\ep(x) J_3\left(x,{z\over h_\ep(x)}\right).$$
As above, we denote $\mathbf{J}_\ep(x,z)=\mathbf{\tilde J}(x,X,Z)$ with
$$\mathbf{\tilde{J}}=(\mathbf{J}, h J_3)\in\RR^{2}\times\RR\,,\,\,\,\,\,\,\mathbf{J}=J_1\mathbf{i}+J_2\mathbf{j}\,.$$
Now we rewrite the governing equations taking into account the above change of variables. This computation is tedious but straightforward so we refer the reader to \cite{Ja7} for technical details.
As a result, the momentum equation \eqref{1,1-J} takes the following form in an $\ep$-independent domain $\Omega\times\RR^2=\{(x,X,Z)\in\RR^2\times\RR^2\times\RR\,:\,x\in (0,1)^2,\,0<Z<1\}$:
\begin{equation}\label{velhor1}
\begin{array}{l}
 \dis -h^2\Delta_x(\mathbf{v}+\mathbf{J}) -{2\over \ep^2}h^2\nabla_x\cdot \nabla_X (\mathbf{v}+\mathbf{J})  -{1\over \ep^4}h^2\Delta_X(\mathbf{v}+\mathbf{J})\\
\noame\dis \qquad\qquad +{2 \over \ep^2}h\,\nabla h \cdot Z\nabla_X\partial_Z (\mathbf{v}+\mathbf{J})+ h\, \Delta h Z\partial_Z (\mathbf{v}+\mathbf{J})  - |\nabla h|^2Z\partial_Z (\mathbf{v}+\mathbf{J}) \\
\noame\dis\qquad\qquad +2 h\,\nabla h\cdot Z\nabla_x\partial_Z (\mathbf{v}+\mathbf{J})- |\nabla h|^2 Z^2\partial_Z^2(\mathbf{v}+\mathbf{J})  - \partial_Z^2(\mathbf{v}+\mathbf{J}) \\
\noame\dis\qquad\qquad+h^2\nabla_x p+ {1\over \ep^2}h^2\nabla_Xp-h\nabla h Z\partial_Z p\\
\noame\dis \qquad\qquad= \dis {2N^2\over \ep^2}h^2\mbox{rot}_X w_3 + 2N^2 h\,\Big(-\partial_Z w_2\, \mathbf{i} + \partial_Z w_1\, \mathbf{j}\Big) + 2N^2 h^2\mbox{rot}_x w_3\,, \end{array}\end{equation}

$$\begin{array}{l}
 \dis -h^2 \Delta_x (v_3+h J_3) -{2 \over \ep^2}h^2\nabla_x\cdot \nabla_X(v_3+h J_3) -{1 \over \ep^4}h^2\Delta_X(v_3+h J_3) +{2 \over \ep^2}h\,\nabla h \cdot Z\nabla_X\partial_Z (v_3+h J_3)\\
\noame\dis \qquad\qquad + h\,\Delta h Z\partial_Z (v_3+h J_3)  - |\nabla h|^2 Z\partial_Z (v_3+h J_3) +2  h\,\nabla h\cdot Z\nabla_x\partial_Z (v_3+h J_3)\\
\noame\dis\qquad\qquad - |\nabla h|^2 Z^2\partial_Z^2(v_3+h J_3)  -\partial_Z^2(v_3+h J_3) + h\,\partial_Z p\\
\noame\dis\qquad\qquad =\dis {2N^2\over \ep^2}h^2\Big(\partial_{X_1}w_2-\partial_{X_2}w_1\Big)
+ 2N^2 h^2 \Big(\partial_{x_1}w_2-\partial_{x_2}w_1 \Big)\,.
 \end{array}$$

The divergence equation (\ref{21,1-J}) is given by
\begin{equation}\label{diverg}h\,\mbox{div}_x \mathbf{v} + {1\over \ep^2}h\mbox{div}_X \mathbf{v}-\nabla h\cdot Z\partial_Z\mathbf{v}+\partial_Z v_3=0\,.
\end{equation}
Finally, the angular momentum equation (\ref{31,1-J}) reads
\begin{equation}\label{microhor}\begin{array}{l}
 \dis
R_c\Big[-\ep^2h^2\Delta_x\mathbf{w} -{2}h^2\nabla_x\cdot \nabla_X \mathbf{w}  -{1\over \ep^2}h^2\Delta_X\mathbf{w} +{2 }h\,\nabla h \cdot Z\nabla_X\partial_Z \mathbf{w} \\
\noame\dis\  + \ep^2 h\, \Delta h Z\partial_Z \mathbf{w}  - \ep^2|\nabla h|^2Z\partial_Z \mathbf{w}  +2 \ep^2 h\,\nabla h\cdot Z\nabla_x\partial_Z \mathbf{w} -\ep^2 |\nabla h|^2 Z^2\partial_Z^2\mathbf{w}  -\ep^2\partial_Z^2\mathbf{w}\Big]+4N^2 h^2\mathbf{w}\\
\noame\dis
\   ={2N^2\over \ep^2}h^2\mbox{rot}_X v_3 + 2N^2 h\,\Big(-\partial_Z (v_2+J_2)\, \mathbf{i} + \partial_Z (v_1+J_1)\, \mathbf{j}\Big) + 2N^2 h^2\mbox{rot}_x v_3\,,
\end{array}\end{equation}
$$\begin{array}{l}
 \dis
R_c \Big[-\ep^2h^2\Delta_x w_3 -{2}h^2\nabla_x\cdot \nabla_X w_3  -{1\over \ep^2}h^2\Delta_Xw_3 +{2}h\,\nabla h \cdot Z\nabla_X\partial_Z w_3 \\
\noame\dis\
 + \ep^2 h\, \Delta h Z\partial_Z w_3  - \ep^2 |\nabla h|^2Z\partial_Z w_3  +2\ep^2 h\,\nabla h\cdot Z\nabla_x\partial_Z w_3-\ep^2|\nabla h|^2 Z^2\partial_Z^2w_3  -\partial_Z^2w_3\Big]+ 4N^2 h^2w_3\\
\noame\dis
\ = {2N^2\over \ep^2}h^2\Big(\partial_{X_1}v_2-\partial_{X_2}v_1\Big)
+ 2N^2 h^2 \Big(\partial_{x_1}v_2-\partial_{x_2}v_1 \Big)\,.
\end{array}$$

%
The boundary conditions (\ref{bc-1})-(\ref{bc-3}) after performing the change of variables remain unchanged:
\begin{eqnarray}\label{bc_zero1}
&&\mathbf{v}=v_3=0,\quad \mathbf{w}=w_3=0\quad\hbox{ on }Z=1,\\
&&v_3=w_3=0\quad\hbox{ on }Z=0,
\end{eqnarray}
whereas conditions (\ref{bc-4})-(\ref{bc-5}) on $Z=0$ read as follows:
\begin{equation}\label{bc_non_resc1}
\begin{array}{rl}
\displaystyle h\,w_2\mathbf{i}-h\,w_1\mathbf{j}=&\displaystyle {\alpha\over 2}\left({h\over \ep^2}\partial_{X_2}(v_3+h J_3)+\partial_Z (v_1+J_1)+h\,\partial_{x_2}(v_3+h J_3)\right)\mathbf{i}\\
\noame&\dis +{\alpha\over 2}\left(-{h\over \ep^2}\partial_{X_2}(v_3+h J_3)+\partial_Z (v_2+J_2)-h\,\partial_{x_2}(v_3+h J_3)\right)\mathbf{j}\,,
\end{array}
\end{equation}

\begin{equation}\label{bc_non_resc2}
\begin{array}{r}\dis
R_c\left(h\,\partial_{X_2}w_3+{\ep}\partial_Z w_1+h\,\ep^2\partial_{x_2}w_3\right)\mathbf{i}
+R_c\left(-h\,\partial_{X_2}w_3+{\ep}\partial_Z w_2-h\,\ep^2\partial_{x_2}w_3\right)
\mathbf{j}\\
\noame\dis =2N^2h\beta(v_2+f_2-s_2)\mathbf{i}-
2N^2h\beta(v_1+f_1-s_1)\mathbf{i}\,.
\end{array}
\end{equation}

Due to the periodic nature of $h_2$, we take into account that
$\tilde{\bf u}$, $\tilde{\bf w}$, $\tilde{\bf J}$ and $p$ are $\mathbb{T}^2$-periodic functions in the variable $X$, namely
\begin{equation}\begin{array}{c}\label{1-periodic}
\tilde{\bf u}(x,X+1,Z)=\tilde{\bf u}(x,X,Z),\quad \tilde{\bf w}(x,X+1,Z)=\tilde{\bf w}(x,X,Z),\\
\noame\dis  \tilde{\bf J}(x,X+1,Z)=\tilde{\bf J}(x,X,Z),\quad p(x,X+1,Z)=p(x,X,Z).
\end{array}
\end{equation}
\subsection{Asymptotic expansion}
We expand the unknowns appearing in the rescaled problem (\ref{velhor1})-(\ref{bc_non_resc2}) as follows:
 \begin{eqnarray}
&&{\mathbf{v}}(x,X,Z)=\mathbf{v}^0(x,X,Z)+ \ep\mathbf{v}^1(x,X,Z)+\ep^2\mathbf{v}^2(x,X,Z)+\cdots\,,\label{exp1}\\[0.2cm]
&&v_3(x,X,Z)=v_3^0(x,X,Z)+ \ep v_3^1(x,X,Z)+\ep^2 v_3^2(x,X,Z)+\cdots\,,\label{exp12}\\[0.2cm]
&&p(x,X,Z)={1\over \ep^2}p^0(x,X,Z)+ {1\over\ep} p^1(x,X,Z)+p^2(x,X,Z)+\cdots\,,\label{exp2}\\[0.2cm]
&&\mathbf{w}(x,X,Z)={1\over \ep} \mathbf{w}^0(x,X,Z)+ \mathbf{w}^1(x,X,Z)+\ep\mathbf{w}^2(x,X,Z)+\cdots\,,\label{exp3}\\[0.2cm]
&&w_3(x,X,Z)={1\over \ep} w_3^0(x,X,Z)+  w_3^1(x,X,Z)+\ep w_3^2(x,X,Z)+\cdots\,.\label{exp31}\end{eqnarray}
Before proceeding, let us determine the asymptotic behavior of the terms involving function $h$. Taking into account that
$$\nabla h(x,X)=\ep\nabla_x h_1(x)+\nabla_X h_2(X),\,\,\,\,\,\,\,\Delta h(x,X)=\ep\Delta_xh_1(x)+{1\over\ep^2}\Delta_X h_2(X),$$
we get
$$\begin{array}{l}
 \dis h^2=\ep^2 h_1^2 + 2\ep^3 h_1 h_2 + \ep^4 h_2^2 \sim O(\ep^2)\,,\\
 \noame\dis
 h\nabla h=\ep h_1\nabla_X h_2+\ep^2 h_1\nabla_x h_1+ \ep^2 h_2\nabla_X h_2 + \ep^3 h_2\nabla_x h_1\sim O(\ep)\,,\\
 \noame\dis
 h\Delta h={1\over \ep}h_1\Delta_X h_2 + h_2\Delta_X h_2+\ep^2 h_1\Delta_x h_1+\ep^3 h_2\Delta_x h_1\sim O\left({1\over \ep}\right)\,,\\
 \noame\dis
 |\nabla h|^2=|\nabla_X h_2|^2+2\ep\nabla_x h_1\nabla_X h_2+\ep^2|\nabla_x h_1|^2\sim O(1)\,.
 \end{array}$$

\subsection{Main order term}
We plug the expansions (\ref{exp1})-(\ref{exp31}) into momentum and divergence equation (\ref{velhor1})-(\ref{diverg}). The main order terms provide
\begin{equation}\label{pr1}\begin{array}{rr}
\dis {1\over \ep^2}:\quad& -h_1^2\Delta_X( \mathbf{v}^0+ \mathbf{J})+h_1^2 \nabla_Xp^0=0\,, \\
\noame\dis
 {1\over \ep^2}:\quad& -h_1^2\Delta_Xv_3^0=0\,,\\
\noame\dis
 {1\over \ep}:\quad& h_1\,\mbox{div}_X\mathbf{v}^0=0\,.
\end{array}\end{equation}
Taking into account the boundary conditions with respect to $X$, we immediately deduce
\begin{equation}\label{expepv0p0v30}
\nabla_X (\mathbf{v}^0+ \mathbf{J})=0\,,\quad\quad \nabla_Xp^0=0\,,\quad\quad  \nabla_X v_3^0=0\,.\end{equation}
In view of that, the main order term from the angular momentum equation (\ref{microhor}) gives
\begin{equation}\label{expep2ep2micro}\begin{array}{rr}
\dis {1\over \ep}:\quad&-R_c h_1^2\Delta_X \mathbf{w}^0 =0\,,\\
\noame\dis
 {1\over \ep}:\quad& -R_c h_1^2\Delta_X w_3^0=0
\end{array}\end{equation}
leading to
 \begin{equation}\label{windp}
\nabla_X\mathbf{w}^0=0\,,\quad\quad\quad \nabla_Xw_3^0=0\,.\end{equation}
We conclude that, at this stage, the problems for velocity/pressure and microrotation are decoupled. As a consequence, we obtain the fast variable independence of the leading order terms.
\subsection{Lower order terms}
The equations satisfied by the lower-order terms from (\ref{velhor1})-(\ref{diverg}) are given by
\begin{equation}\label{corrvel1}\begin{array}{rr}
\dis {1\over \ep}:\quad&  -h_1^2\Delta_X \mathbf{v}^1+h_1^2\nabla_Xp^1+h_1\Delta_X h_2 Z\partial_Z( \mathbf{v}^0+ \mathbf{J})-h_1\nabla_Xh_2 Z\partial_Z p^0=0\,,\\
\noame\dis
 {1\over \ep}:\quad& -h_1^2\Delta_X (v_3^1+h_1 J_3)+h_1\Delta_X h_2 Z\partial_Z v_3^0+h_1\partial_Z p^0=0\,,\\
\noame\dis
 1:\quad& h_1\,\mbox{div}_X \mathbf{v}^1-\nabla_X h_2 \cdot Z\partial_Z (\mathbf{v}^0+ \mathbf{J})+\partial_Z v_3^0=0\,.
\end{array}\end{equation}
We proceed similarly as in \cite{Ja7}. We first consider the last equation (\ref{corrvel1})$_3$ and take the mean value with respect to $X$. Since $h_1$, ${\bf v}^0$, $v_3^0$ are $X$-independent, we deduce $\partial_Z v_3^0=0$ implying  ($v_3^0|_{Z=0,1}=0$)
$$v_3^0=0.$$
Then from (\ref{corrvel1})$_3$ we have
\begin{equation}\label{divv1}
h_1{\rm div}_X({\bf v}^1)=\nabla_X h_2\cdot Z\partial_Z({\bf v}^0+\mathbf{J}).
\end{equation}
Computing the mean value with respect to $X$ in (\ref{corrvel1})$_2$, we arrive at
\begin{equation}\label{v31_X}\partial_Z p^0=0\,,\quad\quad\quad \nabla_X (v_3^1+h_1 J_3)=0.\end{equation}
Finally, applying the curl$_X$ operator of the horizontal component in equation (\ref{corrvel1})$_1$, we obtain
$$h_1{\rm curl}_X\Delta_X ({\bf v}^1)=\nabla_X^{\perp}\Delta_X h_2\cdot Z\partial_Z( {\bf v}^0+\mathbf{J}^0).$$
However, $h_1$ and ${\bf v}^0+\mathbf{J}$ are $X$-independent leading to
\begin{equation}\label{curlv1}
h_1{\rm curl}_X( {\bf v}^1)=\nabla_X^{\perp}h_2\cdot Z\partial_Z ({\bf v}^0+\mathbf{J}).
\end{equation}
Using the identities $\Delta_X ({\bf v}^1)=\nabla_X{\rm div}_X({\bf v}^1)-\nabla_X^{\perp}{\rm curl}_X ({\bf v}^1)$, $\nabla_X\nabla_X h_2-\nabla_X^{\perp}\nabla_X^{\perp}h_2=(\Delta_X h_2){\rm Id}$, we get (in view of (\ref{divv1}), (\ref{curlv1}))
\begin{equation}\label{relationv1v0}
h_1^2\Delta_X( \mathbf{v}^1)=h_1\Delta_X h_2 Z\partial_Z (\mathbf{v}^0+\mathbf{J})\,.
 \end{equation}
Thus, from equation (\ref{corrvel1})$_1$ it remains
$$\nabla_Xp^1=0\,.$$
The next terms from the equations (\ref{velhor1})-(\ref{diverg}) yield
\begin{equation}\label{exp11ep}\begin{array}{rl}
\dis 1:& h_2\Delta_X h_2 Z\partial_Z (\mathbf{v}^0\!+\! \mathbf{J})\!-\!|\nabla_X h_2|^2Z\partial_Z (\mathbf{v}^0+ \mathbf{J})\!-\!|\nabla_X h_2|^2 Z^2\partial_Z^2 (\mathbf{v}^0+ \mathbf{J})\!-\!\partial_Z^2(\mathbf{v}^0+ \mathbf{J})\\
\noame & \dis -2 h_1 h_2\Delta_X(\mathbf{v}^1)+2 h_1\nabla_Xh_2\cdot Z\nabla_X\partial_Z( \mathbf{v}^1)+  h_1\Delta_X h_2 Z\partial_Z(\mathbf{v}^1)) - h_1^2\Delta_X\mathbf{v}^2 \\
\noame&\dis+h_1^2\nabla_xp^0+h_1^2\nabla_Xp^2=2N^2 h_1\left(-\partial_Z w_2^0{\bf i}+\partial_Z w_1^0{\bf j}\right)\,,\\
\noame\dis 1:& -h_1^2\Delta_X (v_3^2) +h_1\Delta_X h_2 Z\partial_Z (v_3^1+ h_1 J_3)+h_1\partial_Z p^1
 =0\,,\\
\noame\dis \ep:& h_1\,\mbox{div}_x(\mathbf{v}^0+ \mathbf{J})+h_1\,\mbox{div}_X (\mathbf{v}^2)+h_2\,\mbox{div}_X(\mathbf{v}^1)-\nabla_x h_1\cdot Z\partial_Z(\mathbf{v}^0+ \mathbf{J})\\
\noame&\qquad\qquad\dis -\nabla_X h_2\cdot Z\partial_Z(\mathbf{v}^1)+\partial_Z (v_3^1+h_1 J_3)=0\,.
\end{array}
\end{equation}
We proceed similarly as above. We take (\ref{exp11ep})$_3$ and compute the mean value in $X$. Using (\ref{divv1}), we get
\begin{equation}\label{Rey11}\mbox{div}_x\left(h_1 (\mathbf{v}^0+\mathbf{J})\right)+\partial_Z\left((v_3^1+h_1 J_3)-\nabla_xh_1\cdot Z (\mathbf{v}^0+\mathbf{J})\right)=0\,.\end{equation}
Now we take the mean value in (\ref{exp11ep})$_2$. In view of (\ref{v31_X})$_2$, we deduce
\begin{equation}\label{Rey112}\partial_Z p^1=0\,.\end{equation}
Last but not least, we compute the mean value with respect to $X$ in (\ref{exp11ep})$_1$. In the following, we address in detail each term of the obtained equation:

$(i)$ Terms involving $\mathbf{v}^0+{\bf J}$: since $\mathbf{v}^0+{\bf J}$ do not depend on $X$, we have\\
 $$ \begin{array}{l}
   \dis \quad \int_{\mathbb{T}^2}\Big[h_2\Delta_X h_2 Z\partial_Z(\mathbf{v}^0+{\bf J})-|\nabla_X h_2|^2 Z\partial_Z (\mathbf{v}^0+{\bf J})
-|\nabla_X h_2|^2 Z^2\partial_Z^2(\mathbf{v}^0+{\bf J})-\partial_Z^2 (\mathbf{v}^0+{\bf J})\Big]dX\\
\noame \dis
\end{array}$$
$$\begin{array}{l} \dis= \quad \left(\int_{\mathbb{T}^2}h_2\Delta_X h_2 \, dX\right)Z\partial_Z(\mathbf{v}^0+{\bf J})-\left(\int_{\mathbb{T}^2}|\nabla_X h_2|^2dX\right) Z\partial_Z (\mathbf{v}^0+{\bf J})
\\
\noame\dis\qquad\qquad-\left(\int_{\mathbb{T}^2}|\nabla_X h_2|^2dX\right) Z^2\partial_Z^2(\mathbf{v}^0+{\bf J})-\partial_Z^2 (\mathbf{v}^0+{\bf J})\\
\noame \dis
\end{array}$$
$$ \begin{array}{l}
 \dis=-\left(\int_{\mathbb{T}^2}|\nabla_Xh_2|^2dX \right)Z\partial_Z(\mathbf{v}^0+{\bf J})-\left(\int_{\mathbb{T}^2}|\nabla_Xh_2|^2dX\right) Z\partial_Z(\mathbf{v}^0+{\bf J})\\
 \noame\dis
 \qquad\qquad\qquad\qquad
 -\left(\int_{\mathbb{T}^2}|\nabla_Xh_2|^2dX\right) Z^2\partial^2_Z(\mathbf{v}^0+{\bf J})-\partial_Z^2(\mathbf{v}^0+{\bf J})\\
 \noame\dis
= -2MZ\partial_Z(\mathbf{v}^0+{\bf J}) -MZ^2\partial_Z^2(\mathbf{v}^0+{\bf J})-\partial_Z^2(\mathbf{v}^0+{\bf J}).
\end{array}$$

It is important to observe that in the above term a new coefficient, denoted by $M$ and defined as
\begin{equation}\label{coefM}
M=\int_{\mathbb{T}^2}|\nabla_X h_2|^2dX
\end{equation}
appears for the first time. This coefficient depends exclusively on the considered rugosity profile and, thus, carries out the effects we seek for in the effective model.\\

$(ii)$ Using the decomposition $\Delta_X{\bf v}^1=\nabla_X{\rm div}_X {\bf v}^1-\nabla_X^\perp{\rm rot}_X{\bf v}^1$, relations (\ref{divv1}), (\ref{curlv1}) and the fact that $(\mathbf{v}^0+{\bf J})$ is $X$-independent, we obtain\\
$$\begin{array}{l}
 \dis T=\int_{\mathbb{T}^2}-2 h_1h_2\Delta_X \mathbf{v}^1
 dX=-\int_{\mathbb{T}^2}2 h_1 h_2\nabla_X{\rm div}_X{\bf v}^1\,dX+\int_{\mathbb{T}^2}2h_1h_2\nabla_X^\perp{\rm rot}_X{\bf v}^1\,dX\\
 \noame\dis\qquad
 =\int_{\mathbb{T}^2}2 h_1 {\rm div}_X{\bf v}^1\nabla_X h_2\,dX+\int_{\mathbb{T}^2}2h_1{\rm rot}_X{\bf v}^1\nabla_X^\perp h_2\,dX
 \\
 \noame\dis\qquad

 =\int_{\mathbb{T}^2}2\left((\nabla_X h_2\cdot Z\partial_Z ({\bf v}^0+{\bf J}))\nabla_X h_2-(\nabla_X^\perp h_2\cdot Z\partial_Z ({\bf v}^0+{\bf J}))\nabla_X h_2^\perp\right)dX \\
 \noame
 \dis
 \qquad
  =\int_{\mathbb{T}^2}2|\nabla_X h_2|^2 Z\partial_Z({\bf v}^0+{\bf J})\,dX
 \\
 \noame\dis
 \qquad=2M Z\partial_Z ({\bf v}^0+{\bf J})\,.
 \end{array}$$

$(iii)$ The remaining terms involving ${\bf v}^1$: arguing as for the previous term, we get
$$\begin{array}{l}\dis \int_{\mathbb{T}^2}\Big[2h_1\nabla_Xh_2\cdot  Z\nabla_X\partial_Z\mathbf{v}^1+h_1\Delta_X h_2 Z\partial_Z \mathbf{v}^1\Big]dX\\
= \dis \int_{\mathbb{T}^2}\Big[-2h_1h_2\cdot  Z\partial_Z(\Delta_X\mathbf{v}^1)+h_1 h_2 Z\partial_Z (\Delta_X\mathbf{v}^1)\Big]dX\\
\noame  \dis
%
= -\int_{\mathbb{T}^2}h_1 h_2Z\partial_Z\Delta_X\mathbf{v}^1 
dX\\
\noame\dis
=Z\partial_Z T= MZ\partial_Z({\bf v}^0+{\bf J}) + MZ^2\partial_Z^2({\bf v}^0+{\bf J}).
\end{array}$$

(iv) Terms involving ${\bf v}^2$ and $p^2$: integrating by parts, because the function  $h_1$ depends only on $x$ and the unknowns ${\bf v}^2$, $p^2$ are periodic in $X$, it follows
$$-\int_{\mathbb{T}^2}h_1^2\Delta_X \mathbf{v}^2dX+\int_{\mathbb{T}^2}h_1^2\nabla_Xp^2dX
=0\,.$$
Adding all contributions $(i)-(iv)$, and recalling that $p^0$, ${\bf w}^0$ do not depend on $X$, from (\ref{exp11ep})$_1$ it is straightforward to deduce
\begin{equation}\label{limit-eq-vel}-\partial_Z^2(\mathbf{v}^0+{\bf J})+MZ\partial_Z(\mathbf{v}^0+{\bf J}) + h_1^2\nabla_x p^0=2N^2h_1\left(-\partial_Z w_2^0\,{\bf i} + \partial_Zw_1^0\,{\bf j}\right)\,.
\end{equation}

Now, we turn our attention to microrotation equation (\ref{microhor}) and compute the lower-order terms. In view of (\ref{windp}), we get
\begin{equation}\label{expepepmicro}\begin{array}{rl}
\dis {1}:\quad&\dis    R_c\Delta_X \mathbf{w}^1=R_c{\Delta_Xh_2\over h_1}Z\partial_Z\mathbf{w}^0\,,\\
\noame\dis
 {1}:\quad&\dis  R_c\Delta_X w_3^1=R_c{\Delta_X h_2\over h_1} Z\partial_Zw_3^0\,,
\end{array}\end{equation}
\begin{equation}\label{exp11micro}\begin{array}{rl}
\dis {\ep}:\quad&    R_c\Big[h_2\Delta_X h_2 Z\partial_Z\mathbf{w}^0-|\nabla_X h_2|^2 Z\partial_Z \mathbf{w}^0
-|\nabla_X h_2|^2 Z^2\partial_Z^2\mathbf{w}^0-\partial_Z^2 \mathbf{w}^0\Big]\\
\noame &\dis
+R_c\Big[-2 h_1h_2\Delta_X \mathbf{w}^1+2h_1\nabla_Xh_2 Z\nabla_X\partial_Z\mathbf{w}^1+h_1\Delta_X h_2 Z\partial_Z \mathbf{w}^1\Big]\\
\noame & \dis
-R_c h_1^2\Delta_X \mathbf{w}^2+ 4N^2 h_1^2 {\bf w}^0=2N^2 h_1\left(-\partial_Z(v_2^0+J_2){\bf i} + \partial_Z(v_1^0+J_1){\bf j}\right)\,,\\[0.2cm]
\noame\dis
 {\ep}:\quad& R_c\Big[h_2\Delta_X h_2 Z\partial_Z w_3^0-|\nabla_X h_2|^2 Z\partial_Z w_3^0-|\nabla_X h_2|^2 Z^2\partial_Z^2\mathbf{w}^0-\partial_Z^2 w_3^0\Big]\\
 \noame
 \noame &\dis+R_c\Big[-2h_1h_2\Delta_X w_3^1+2h_1\nabla_X h_2\cdot Z\nabla_X\partial_Zw_3^1+h_1\Delta_X h_2 Z\partial_Zw_3^1\Big]=0. 
\end{array}\end{equation}
Taking the mean value in (\ref{exp11micro})$_1$ with respect to $X$, using relation (\ref{expepepmicro})$_1$, and proceeding using similar arguments as above, we arrive at
\begin{equation}\label{limit_eq_mic}-R_c\partial_Z^2\mathbf{w}^0+R_cMZ\partial_Z\mathbf{w}^0+4N^2 h_1^2 {\bf w}^0=2N^2h_1\left(-\partial_Z(v_2^0+J_2){\bf i} + \partial_Z(v_1^0 + J_1){\bf j}\right)\,.\end{equation}
Proceeding analogously for (\ref{exp11micro})$_2$, we derive the equation satisfied by $w_3^0$:
$$-R_c\partial_Z^2w_3^0+R_cMZ\partial_Zw_3^0=0
$$
providing  (recall $w_3^0|_{Z=0,1}=0$)
$$w_3^0=0\,.$$

\begin{rem}
The main order terms from the boundary conditions (\ref{bc_non_resc1})-(\ref{bc_non_resc2}) provide nothing relevant to the formal development. Indeed, since $v_3^0=w_3^0=0$, the relations satisfied by those terms are trivially fulfilled. The lower-order terms from (\ref{bc_non_resc1})-(\ref{bc_non_resc2}) yield:
\begin{equation}\label{bcexp1}\begin{array}{ll}
\dis {1}:\quad&\dis \partial_Z(v_1^0+J_1)\,\mathbf{i} + \partial_Z(v_2^0+J_2)\mathbf{j} ={2\over \alpha} h_1\Big[w_2^0\,{\bf i}+ w_1^0\,{\bf j}\Big],\\
\noame\dis
%
\dis {1}:\quad&\dis R_c\Big[\partial_Zw_1^0\mathbf{i} +\partial_Zw_2^0 \mathbf{j}\Big]= 2N^2h_1\beta\Big[(v_2+f_2-s_2)\,\mathbf{i}-(v_1+f_1-s_1)\,\mathbf{j}\Big].
\end{array}\end{equation}
\end{rem}
\ \\

\subsection{Effective system}
Using (\ref{v-u-J}), from (\ref{Rey11}), (\ref{limit-eq-vel}) and (\ref{limit_eq_mic}) we deduce the system satisfied by the effective velocity, pressure and microrotation:
\begin{equation}\label{effective}\left\{\begin{array}{rl}
-\partial_Z^2\mathbf{u}^0+MZ\partial_Z\mathbf{u}^0+ h_1^2\nabla_x p^0-2N^2h_1\left(-\partial_Z w_2^0\,{\bf i} + \partial_Zw_1^0\,{\bf j}\right)=0\,,\\
\noame\dis \partial_Zp^0=0\,,\\
\noame\dis{\rm div}_x\int_0^1h_1{\bf u}^0 dZ=0\,,\\
\noame\dis
-R_c\partial_Z^2\mathbf{w}^0+R_cMZ\partial_Z\mathbf{w}^0+4N^2 h_1^2 {\bf w}^0-2N^2h_1\left(-\partial_Zu_2{\bf i} + \partial_Z u_1{\bf j}\right)=0\,.
\end{array}\right.\end{equation}
We complete the above system with the boundary conditions (see (\ref{bc_zero1}), (\ref{bcexp1})):
\begin{equation}\label{bceffectivew1}
\mathbf{u}^0=0\ \hbox{ on } Z=1,\quad \mathbf{w}^0=0\ \hbox{ on } Z=1,
\end{equation}
\begin{equation}\label{bceffectivew21}
\partial_Zu_1^0={2\over \alpha} h_1 w_2^0 \ \hbox{ on } Z=0,\quad \partial_Z u_2^0={2\over \alpha}h_1 w_1^0 \ \hbox{ on } Z=0, \\
\end{equation}
\begin{equation}\label{bceffectivew22}
R_c\,\partial_Z w_1^0=2N^2h_1\beta(u_2^0-s_2) \ \hbox{ on } Z=0,\quad R_c\,\partial_Z w_2^0=-2N^2h_1\beta(u_1^0-s_1)\ \hbox{ on } Z=0\,.
\end{equation}
The problem (\ref{effective})-(\ref{bceffectivew22}) forms the effective system describing the macroscopic flow. It is essential to observe that the effective equations (\ref{effective})$_1$, (\ref{effective})$_4$ are explicitly corrected by the roughness-induced term containing the coefficient $M$. Indeed, by putting $M=0$ (i.e.~no roughness introduced) and assuming the invariance in the $x_2$-direction, we obtain the same limit equations as derived in \cite{Bayada4} via weak convergence method. Moreover, comparing our result with the one derived in \cite{Ja7} for simple zero boundary conditions, we clearly detect the effects of the nonzero (physically relevant) boundary conditions imposed on the governing problem.

Let us mention that the asymptotic model (\ref{effective})-(\ref{bceffectivew22}) can also be derived using two-scale convergence method (proposed originally by Allaire \cite{Allaire}). The procedure would follow exactly the same arguments as presented in \cite{Ja7} providing the rigorous justification of our formally obtained model. Instead of that, here we prefer to perform numerical simulations of the above system focusing on the effects of the rough boundary. The main difficulties arise due to the roughness-induced terms in the effective equations preventing us to solve them analytically (as it was possible in \cite{Bayada4}). For that reason, we need to employ a different numerical strategy, as presented in the following section.




\section{Roughness-induced effects: numerical results}


In this section, we use the asymptotic system derived in Sec.~3 to investigate the performance of a linear slider bearing with a rough inclined surface, lubricated with a micropolar fluid. We denote by $(\widehat{x}, \widehat{z})\in \RR^2\times \RR$ the space variable, with $\widehat{x} = (\widehat{x}_1,\widehat{x}_2)$. The lower surface is a plane wall at $\widehat{z}=0$ that moves with constant velocity $g$ in the $\widehat{x}_1$ direction.
We assume that the bearing is invariant in the $\widehat{x}_2$ direction, and denote by $c$ the maximum distance between the surfaces, by $L$ the length of the device and by $m$ the slope of the inclined surface. The geometry of the domain $\widehat{\Omega}_\ep$ is described by the gap function
\[
\widehat h_\ep(\widehat x_1) = c\left( 1+m\frac{\widehat x_1}{L} \right) + \frac{c^2}{L}\ h_2\left( \frac{L}{c^2}\widehat x_1 \right),\,\,\,\,\,\,\,\,\widehat x_1\in (0,L)\,.
\]
Here function $h_2$ is a regular function, periodic with period $1$ and with zero average, describing the roughness of the inclined surface. In view of that, the domain reads
\[
\widehat{\Omega}_\ep = \left\{ (\widehat x, \widehat z)\in \RR^2\times \RR\ :\ \widehat x_1\in(0,L),\ 0<\widehat z < \widehat h_\ep(\widehat x_1) \right\},\,\,\,\,\,\,\,\,\ep = \frac{c}{L}\,.
\]
Introducing non-dimensional space variables (as in Sec.~2)
\[
x=\frac{\widehat x}{L},\quad z=\frac{\widehat z}{L},
\]
the non-dimensional gap function $h_\ep=\frac{\widehat h_\ep}{L}$ is given by
\[
h_\ep(x_1) = \ep h_1(x_1) + \ep^2\ h_2\left(\frac{x_1}{\ep^2}\right),\qquad x_1\in (0,1),
\]
where $h_1(x_1)=1+mx_1$. The (rescaled) domain now has the form
\[
\Omega_\ep=\left\{(x,z)\in \RR^2\times \RR,\ 0<x_1<1,\ 0<z<h_\ep(x_1)\right\}.
\]
Due to the invariance in the $x_2$-direction, the velocity $\ub_\ep$ and the microrotation $\wb_\ep$ reduce to
$$\ueps=(u_{1,\ep},0,u_{3,\ep}),\qquad \weps=(0,w_{2,\ep},0)$$
leading to the following limit system posed in $\Omega:=(0,1)^2$ (see (\ref{effective})-(\ref{bceffectivew22})):
\begin{equation}\label{effective_u_1_w_2_2D}\left\{\begin{array}{rl}
-\partial_Z^2 u_1+MZ\partial_Zu_1+ h_1^2d_{x_1} p+2N^2h_1\partial_Z w_2=0\,,\\
\noame\dis
-R_c\partial_Z^2w_2+R_cMZ\partial_Zw_2+4N^2 h_1^2 w_2-2N^2h_1 \partial_Z u_1=0\,,\\
\noame\dis
d_{x_1}\int_0^1h_1 u_1\,dZ=0\\
\noame\dis
u_1=0\ \hbox{ on } Z=1,\quad w_2=0\ \hbox{ on } Z=1,\\
\noame\dis
\partial_Zu_1={2\over \alpha} h_1 w_2 \ \hbox{ on } Z=0,\quad  R_c\,\partial_Z w_2=-2N^2h_1\beta(u_1-s_1) \ \hbox{ on } Z=0
\end{array}\right.\end{equation}

In order to be able to compare our results with those from \cite{Bayada4}, we may transform $\Omega$ in the following set
$$\Omega_l=\left\{(x_1,Y),\ 0<x_1<1,\ 0<Y<h_1(x_1)\right\},$$
by introducing the change of variables $Z=Y/h_1(x_1)$. This yields the following system posed in $\Omega_l$:
\begin{equation}\label{effective_u_1_w_2_2D_cha}\left\{\begin{array}{rl}
-\partial_Y^2 u_1+{M\over h_1^2}Y\partial_Yu_1+ d_{x_1} p+2N^2\partial_Y w_2=0\,,\\
\noame
-R_c\partial_Y^2w_2+{R_cM\over h_1^2}Y\partial_Yw_2+4N^2  w_2-2N^2 \partial_Z u_1=0\,,\\
\noame\dis
d_{x_1}\int_0^{h_1} u_1\,dY=0,\\
\noame\dis
u_1=0\ \hbox{ on } Y=h_1(x_1),\quad w_2=0\ \hbox{ on } Y=h_1(x_1),\\
\noame\dis
\partial_Yu_1={2\over \alpha}  w_2 \ \hbox{ on } Y=0,\quad  R_c\,\partial_Y w_2=-2N^2\beta(u_1-s_1)  \ \hbox{ on } Y=0\,.
\end{array}\right.
\end{equation}

Note that in the case $M=0$ (i.e.~in absence of rugosity), the asymptotic system \eqref{effective_u_1_w_2_2D_cha} is exactly the model studied numerically in \cite{Bayada4}. The terms ${M\over h_1^2}Y\partial_Yu_1$ and ${R_cM\over h_1^2}Y\partial_Yw_2$ appearing in the momentum equations, are the expression of the roughness effect. Since these terms are not exact differentials, we are not able to solve \eqref{effective_u_1_w_2_2D_cha} (or equivalently, \eqref{effective_u_1_w_2_2D}) analytically. To circumvent this difficulty, we propose to approximate the coupled unknown functions $u_1, w_2$ by an iterative scheme.
\subsection{Approximation of the limit system \eqref{effective_u_1_w_2_2D} by an iterative scheme}\label{Subsection:IterativeScheme}
We denote by $V_Z$ the space
$
V_Z = \left\lbrace v\in H^1(0,1),\ v(1) = 0	\right\rbrace
$
endowed with the norm $\|v\|_{V_Z}:=(\int_0^1 |d_Z v|^2\ dZ)^{1/2}$. Starting from $\uzero$, we construct a sequence of approximations of $u_1,w_2$ as follows. Assume that $\un$ is given in $\LV$. For every $x_1\in (0,1)$, we define $\wn(x_1,\cdot)\in V_z$ as the variational solution of the following problem:
\begin{align}
-R_c\d_Z^2 \wn+R_cMZ\d_Z\wn + 4N^2 h_1^2 \wn=2N^2h_1 \partial_Z \un\,,\label{eq-wn}\\
\wn = 0 \ \hbox{ on } Z=1,\quad R_c\,\partial_Z \wn = -2N^2h_1\beta(\un-s_1) \ \hbox{ on } Z=0. \label{CL-wn}
\end{align}
To construct the new approximation $\unp$, we need to deal with the constraint
\be
d_{x_1}\int_0^1 h_1 \unp\,dZ=0. \label{constraint-unp}
\ee
To this aim, we will use a potential $\psi\in H^1_Z(0,1)$, solution of
\begin{align}
-d^2_{Z}\psi + MZ\ d_Z\psi = 1\qquad \forall Z\in (0,1), \label{eq-psi}\\
d_Z\psi(0) = 0,\quad \psi(1) = 0.\label{CL-psi}
\end{align}
We denote by $\bar{\psi}=\int_0^1 \psi(Z)\ d Z$ its average.
 Let $\unptilde$ be the solution to the following unconstrained problem:
\begin{align}
-\partial_Z^2 \unptilde+MZ\partial_Z\unptilde = -2N^2h_1\partial_Z \wn\,, \label{eq-unptilde}\\
\unptilde=0\ \hbox{ on } Z=1,\quad \partial_Z \unptilde = {2\over \alpha} h_1 \wn \ \hbox{ on } Z=0. \label{CL-unptilde}
\end{align}
We define the pressure $\pnp(x_1)$ as the solution to the following Reynolds equation:
\be
d_{x_1}\left(\psibar\ h_1^3\ d_{x_1}\pnp \right) = d_{x_1}\int_0^1 h_1\ \unptilde\ dZ,
\label{Reynolds-pnp}
\ee
completed with Dirichlet boundary conditions
\be
\pnp(0) = \pnp(1) = 0.
\label{CL-pnp}
\ee
Such boundary conditions for the pressure are usual in the lubrication field.
Then, we set
\begin{equation}
\unp = \unptilde -h_1^2 d_{x_1}\pnp \psi\qquad \forall x_1\in (0,1),\quad \hbox{a.e.}\ Z\in (0,1).
\label{Def:unp}
\end{equation}
By construction, $\unp,\pnp$ satisfy the equation
\be
-\partial_Z^2 \unp+MZ\partial_Z\unp+ h_1^2d_{x_1} \pnp+2N^2h_1\partial_Z \wn=0 \label{eq-unp}
\ee
as well as the boundary conditions
\be
\unp=0\ \hbox{ on } Z=1,\quad \partial_Z \unp = {2\over \alpha} h_1 \wn \ \hbox{ on } Z=0. \label{CL-unp}
\ee
Moreover, it is easy to check that the constraint \eqref{constraint-unp} is fulfilled.

\subsection{Convergence result}

\paragraph{Stability condition} We begin by proving a stability result for the sequences $(\un),(\wn)$ in $\LV$. Let us first recall the Poincar\'e inequality in $V_Z$:
\be
\forall \varphi\in V_Z\quad \int_0^1 \varphi(Z)^2\ dZ \leq \int_0^1 (d_Z\varphi(Z))^2\ dZ.
\label{Ineq-Poincare-VZ}
\ee
Using the identity $(d_Z\varphi) \varphi = \frac{1}{2}d_Z(\varphi^2)$ and an integration by parts, one can deduce from \eqref{Ineq-Poincare-VZ} the following coercivity inequality:
\be \label{Ineq-coercivity}
\forall \varphi\in V_Z\quad \int_0^1(d_Z\varphi)^2 + M\int_0^1 Z\ d_Z\varphi\ \varphi \geq \left(1-\frac{M}{2}\right) \int_0^1(d_Z\varphi)^2.
\ee
Multiplying \eqref{eq-wn} by $\wn$, integrating by parts and taking into account boundary conditions \eqref{CL-wn}, we obtain the relation
\begin{align*}
R_c\int_0^1 (\d_Z\wn)^2\ dZ + R_c M\int_0^1 Z\ \d_Z\wn\ \wn \ dZ + 4N^2h_1^2\int_0^1 (\wn)^2\ dZ \nonumber \\
= 2N^2h_1\int_0^1 \d_Z\un\ \wn\ dZ + 2N^2h_1\beta(\un-s)_{|Z=0}\ {\wn}_{|Z=0}.
\end{align*}
Using the trace inequality
\bestar
\varphi(0)\leq \left(\int_0^1 (\d_Z\varphi)^2\ dZ \right)^{1/2} \quad \forall \varphi\in V_Z,
\eestar
H\"older inequality and inequalities \eqref{Ineq-Poincare-VZ} and \eqref{Ineq-coercivity}, we deduce the following estimate:
\be
\left(
\int_0^1 (\d_Z \wn)^2\ dZ
\right)^{1/2} \leq
\frac{2N^2}{R_c(1-\frac{M}{2})} (\sup h_1) \left(
\left(
\int_0^1 (\d_Z \un)^2\ dZ
\right)^{1/2}
+\beta \left| {\un} - s_1\right|_{|Z=0}
\right)
.
\label{Estimate-wn}
\ee
By the same arguments, we obtain the analogous inequality for $\unptilde$:
\bestar
\left(
\int_0^1 (\d_Z \unptilde)^2\ dZ
\right)^{1/2} \leq
\frac{2N^2+\frac{2}{\alpha}}{1-\frac{M}{2}}
(\sup h_1)
\left(
\int_0^1 (\d_Z \wn)^2\ dZ
\right)^{1/2},
\eestar
and we deduce from \eqref{Estimate-wn}:
\be
\left(
\int_0^1 (\d_Z \unptilde)^2\ dZ
\right)^{1/2} \leq
\frac{2N^2\left(2N^2+\frac{2}{\alpha}\right)}{\left(1-\frac{M}{2}\right)^2}\ (\sup h_1)^2
\left[
\left(
\int_0^1 (\d_Z \un)^2\ dZ
\right)^{1/2}
+\beta \left| {\un} - s_1\right|_{|Z=0}
\right].
\label{Estimate-unptilde}
\ee
To estimate the pressure term $\pnp$, we multiply \eqref{Reynolds-pnp} by $\pnp$ and integrate by parts to obtain
\be
\int_0^1 (d_{x_1}\pnp)^2\ dx_1 \leq 2\left[\frac{\sup h_1}{\psibar(\inf h_1)^3}\right]^2 \int_0^1  \int_0^1 (\d_Z\unptilde)^2\ dZ dx_1.
\label{Estimate-pnp}
\ee
Using \eqref{Def:unp}, \eqref{Estimate-unptilde} and \eqref{Estimate-pnp}, we finally obtain
\be
\| \unp \|_{\LV} \leq C\left[ (1+\beta)\|\un\|_{\LV} + \beta s_1 \right], \label{Estimate-unp}
\ee
where the constant $C=C(\alpha,M,h_1)$ is given by
\[
C = \sqrt{2}\ \frac{2N^2\left(2N^2+\frac{2}{\alpha}\right)}{\left( 1-\frac{M}{2} \right)^2}
(\sup h_1)^2\left[ 1+\sqrt{2}\frac{\|\psi\|_{V_Z}}{\psibar}\left(\frac{\sup h_1}{\inf h_1}\right)^3 \right].
\]
We conclude that if the condition
\beq
C(1+\beta) \leq 1
\label{Condition:stabilite-schema}
\eeq
holds true, then the sequence $(\un)$ is bounded in $\LV$. In that case, estimate \eqref{Estimate-wn} yields that $(\wn)$ is bounded in the same space.

\paragraph{Convergence of the scheme}
Assume that condition \eqref{Condition:stabilite-schema} is true. The space $\LV$ being a Hilbert space, there exist $u_1,w_2\in \LV$ such that, up to a subsequence,
\[
\un \wra u_1,\quad \wn \wra w_2\quad \mbox{weakly in }\LV.
\]
Since $V_Z$ is continuously embedded in $L^2(0,1)$, it is easy to see that the injection $\LV\embed \LL$ is continuous. As a result, the weak convergence still holds in $\LL$:
\[
\un \wra u_1,\quad \wn \wra w_2\quad \mbox{weakly in }\LL.
\]
Applying the same argument to the trace operator from $V_Z$ on $Z=0$, we obtain
\[
\un(\cdot,0) \wra u_1(\cdot,0),\quad \wn(\cdot,0) \wra w_2(\cdot,0)\quad \mbox{weakly in }L^2_{x_1}(0,1).
\]

\paragraph{Obtention of the limit equations}
We multiply equation  \eqref{eq-wn} by a function $\xi\in\LV$. Integrating by parts in $Z\in (0,1)$, using the boundary conditions \eqref{CL-wn} and integrating in $x_1\in (0,1)$, we get
\begin{align*}
& R_c\int_0^1\int_0^1 \d_Z\wn\ \d_Z\xi\ dZdx_1 + R_c M\int_0^1\int_0^1 Z\d_Z\wn\ \xi \ dZdx_1\\
&+4N^2\int_0^1 h_1^2\int_0^1 \wn\ \xi\ dZdx_1\\
& = 2N^2\left\lbrace \int_0^1\int_0^1 h_1\ \d_Z\un\ \xi\ dZdx_1\right.
 \left.+ \int_0^1h_1(x_1)\ \beta(\un(x_1,0)-s)\ \xi(x_1,0)\ dx_1\right\rbrace.
\end{align*}
Notice that the terms $h_1$ or $h_1^2$ that appear in certain integrals do not play any role in the convergence of the corresponding terms. As a consequence, one can pass to the limit in every term of this weak formulation to conclude that
\begin{align*}
& R_c\int_0^1\int_0^1 \d_Zw_2\ \d_Z\xi\ dZdx_1 + R_c M\int_0^1\int_0^1 Z\d_Zw_2\ \xi \ dZdx_1\\
&+4N^2\int_0^1 h_1^2\int_0^1 w_2\ \xi\ dZdx_1\\
& = 2N^2\left\lbrace \int_0^1\int_0^1 h_1\ \d_Zu_1\ \xi\ dZdx_1\right.
 \left.+ \int_0^1h_1(x_1)\ \beta(u_1(x_1,0)-s)\ \xi(x_1,0)\ dx_1\right\rbrace.
\end{align*}
The passing to the limit in the equation for $\unp$ is similar, except that we have to deal with the pressure term $d_{x_1}p^{n+1}$. To this aim, we combine estimates \eqref{Estimate-pnp} and \eqref{Estimate-unptilde} to obtain that $d_{x_1}p^{n+1}$ is bounded in $L^2_{x_1}(0,1)$. Thus, by continuity of the injection $L^2_{x_1}(0,1)\embed \LL$, there exists $p\in L^2_{x_1}(0,1)$ such that, up to a subsequence,
\[
d_{x_1}p^{n+1}\wra d_{x_1}p\quad \mbox{weakly in }L^2_{x_1}(0,1) \mbox{ and in }\LL.
\]
This allows to obtain the relation:
\begin{align*}
& \int_0^1\int_0^1 \d_Z u_1\ \d_Z\xi\ dZdx_1 + M\int_0^1\int_0^1 Z\d_Z u_1\  \xi \ dZdx_1\\
& + \int_0^1\int_0^1 h_1^2\ d_{x_1}p\ \xi\ dZdx_1 \\
& = -2N^2\left\lbrace \int_0^1\int_0^1 h_1\ \d_Z w_2\ \xi\ dZdx_1\right.
 \left.- \int_0^1h_1(x_1)\ \frac{2}{\alpha}w_2(x_1,0)\ \xi(x_1,0)\ dx_1\right\rbrace.
\end{align*}
Finally, to deal with the constraint, we multiply the relation \eqref{constraint-unp} by function $\zeta\in H^1_0(0,1)$ and integrate by part. This yields
$
\int_0^1 \int_0^1 h_1 \unp\ d_{x_1}\zeta\ dZdx_1 =0
$. By weak convergence of $\unp$ in $\LL$, we deduce that
\[
\int_0^1 \int_0^1 h_1 u_1\ d_{x_1}\zeta\ dZdx_1 =0.
\]
Consequently, the constraint $d_{x_1}\int_0^1 h_1 u_1\,dZ=0$ is satisfied in $H^{-1}(0,1)$.

\subsection{Numerical results}
In order to study numerically the influence of the roughness on the performance of the slider bearing, we simulate the asymptotic system \eqref{effective_u_1_w_2_2D} using different values of the parameter $M$ that describes the main order effect of the rugosity. As already mentioned, the case $M=0$ was investigated in \cite{Bayada4}, where the corresponding system was solved analytically. So as to facilitate the comparison with the results obtained by these authors, we set the maximum distance between the plates to $c=1$ and the slope of the inclined surface to $m=-0.5$, so the function $h_1$ reads
\[
h_1(x_1) = 1-0.5x_1,\quad x_1\in (0,1).
\]
We define the characteristic velocity of the plane surface by $g_0=g$, where $g$ is the imposed velocity of the wall; this yields $s_1=1$. Also,
we introduce the same parameters
\[
\bar{\nu}_b = \frac{1-\alpha N^2}{1-N^2},\quad \delta = \frac{R_c}{2N^2\beta}.
\]
With this convention, the solution to system \eqref{effective_u_1_w_2_2D} depends on $5$ dimensionless parameters: $N, R_c, \bar{\nu}_b,\delta$ and the roughness parameter $M$.

In the mechanical engineering literature, the following quantities appear to be relevant to characterize the efficiency of linear bearings:
\begin{itemize}
\item the pressure profile $p$,
\item the load carrying capacity $W$,
\item the coefficient of friction $c_f$ along the moving surface.
\end{itemize}
In the context of the lubrication limit developed in Sec.~3, $W$ and $c_f$ can be expressed in a dimensionless setting by using the solution $u_1,w_2,p$ to the effective system \eqref{effective_u_1_w_2_2D}. The load reads
\[
W=\int_0^1 p(x_1)\,dx_1.
\]
Introducing the dimensionless friction force
\[
F=\int_0^1\left({\partial u_1\over \partial Y}(x_1,0)-2N^2 w(x_1,0)\right)\,dx_1,
\]
the coefficient of friction is then defined by
\[
c_f = \frac{F}{W}.
\]
We denote by $W_0,F_0, {c_f}_0$ the corresponding values, in the case $M=0$, that we obtain applying the analytical formulas derived in \cite{Bayada4}.

\paragraph{Implementation of the iterative scheme}

To simulate system \eqref{effective_u_1_w_2_2D}, we have implemented the iterative scheme presented in Sec.~\ref{Subsection:IterativeScheme}. To this aim, we have used a finite-element approach to approximate the different variational problems \eqref{eq-wn}-\eqref{CL-wn}, \eqref{eq-psi}-\eqref{CL-psi}, \eqref{eq-unptilde}-\eqref{CL-unptilde} and \eqref{Reynolds-pnp}-\eqref{CL-pnp}.

Define a step $h_Z=\frac{1}{n_Z+1}$ and a discretization $0=Z_0<\ldots<Z_{n_Z+1}=1$ of the (vertical) segment $[0,1]$ with $Z_k = kh_Z$, $0\leq k\leq n_Z+1$. We introduce the space
\[
V_{n_Z} = \left\{ v\in C([0,1]),\quad v(1)=0,\quad \forall k\in[|0,n_Z|]\quad v_{|[Z_k,Z_{k+1}]}\ \textrm{is affine.}  \right\}
\]
For a fixed value of $M$, we approximate the potential $\psi$ by computing $\psi_{app}\in V_{n_Z}$, solution of the variational formulation
\[
\forall \varphi\in V_{n_Z}\quad \int_0^1 \left(d_Z\psi_{app}\ d_Z\varphi + M\ Z\ d_Z \psi_{app}\ \varphi\right)\ d_{x_1} = \int_0^1 \varphi\ dx_1.
\]
The treatment of problems \eqref{eq-wn}-\eqref{CL-wn} and \eqref{eq-unp}-\eqref{CL-unp} requires to address problem \eqref{Reynolds-pnp}-\eqref{CL-pnp} simultaneously. To this aim, we introduce in a similar way  a regular subdivision $0=x_1^0<x_1^1<\ldots < x_1^{n_1}=1$ of the (horizontal) segment $[0,1]$, defined by the step $h_{x_1}=\frac{1}{n_1+1}$ and the formula $x_1^i = ih_{x_1}$, $0\leq i \leq n_1+1$. At each step $n$, the approximation of $\pnp$ will be sought in the space
\[
M_{n_1} = \left\{ q\in C([0,1]),\quad q(0)=q(1)=0,\quad \forall i\in[|0,n_1|]\quad q_{|[x_i,x_{i+1}]}\ \textrm{is affine.}  \right\}
\]
Since the right-hand side of equation \eqref{Reynolds-pnp} is defined by the derivative of
$\int_0^1 h_1\ \unptilde\ dZ$ with respect to the horizontal coordinate $x_1$, we locate the approximations of $\unp$ (and consequently, $\wnp$) along the segments $\{x_{i+\frac{1}{2}}\}\times[0,1]$, where $x_{i+\frac{1}{2}}=(i+\frac{1}{2})h_{x_1}$. For each middle point $x_{i+\frac{1}{2}}$, we solve problems \eqref{eq-wn}-\eqref{CL-wn} and \eqref{eq-unp}-\eqref{CL-unp} by a finite element method, similar to one used to compute the potential $\psi$. In particular, the integrals  $\int_0^1 h_1\ \unptilde\ dZ$ can be computed then at points $x_{i+\frac{1}{2}}$, and their derivatives can be approximated at grid points $x_i$ consistently, using a classical centered scheme.

\paragraph{Numerical results and comments}

In the simulations, we have used the following constant parameters:
\[
N = 0.1,\ R_c = 0.01,\ \nu_b = 0.1,\ \delta = 0.01.
\]
In Fig.~\ref{Fig:p1}, Fig.~\ref{Fig:p2}, Fig.~\ref{Fig:p3}, we have plotted the pressure against $x_1$, with $N=0.1, N=0.2$ and $N=0.3$, respectively. In each case, we have used $M=0$, $M=0.5$, $M=1$. We notice that, for each choice of $N$, the maximum value of $p$ increases when $M$ increases. However, this effect appears to be more pronounced when the value of parameter $N$ is small. As a direct consequence, we observe in Fig.~\ref{Fig:Load} that the relative load $W/W_0$ is an increasing function of $M$. Again, the effect of an increase of $M$ gets more pronounced as $N$ diminishes, for this choice of values of $N$.

Note that $W_0$ was computed using a totally different approach than the one proposed in \cite{Bayada4}; in the case $M=0$, both methods give the same results. This was a major point in the validation of our numerical strategy. This property holds for the computation of the coefficient of friction. In Fig.~\ref{Fig:FrictionCoef}, we have plotted the relative coefficient of friction $c_f/{c_f}_0$ against $M$, with the same choice of $N$. We observe that the increase of $M$ decreases the coefficient of friction; however, the influence of the parameter $N$ is more debatable here.

As a conclusion, we can deduce from these numerical examples that the characteristics of a linear slider bearing with a rough inclined surface  can be modified by the presence of the rugosity, in the following way. With respect to bearing equipped with a plane surface, a bearing with a rough surface may exhibit a larger load carrying capacity, as well as a lower coefficient of friction. Thus, the use of rough surfaces may contribute to enhance the mechanical performance of such device.

\begin{figure}
\begin{center}
\begin{tikzpicture}[scale=0.8]
\begin{axis}[ axis x line=bottom, axis y line = left,
xlabel={$x_1$}, ylabel={pressure $p$},
ymin=0, ymax=0.12,
y tick label style={
        /pgf/number format/.cd,
            fixed,
            fixed zerofill,
            precision=2,
        /tikz/.cd
    }
,
legend entries={$M=0$, $M=0.5$, $M=1$}, legend style={at={(0.02,1)}, anchor=north west}]
\addplot[mark = o] table[x index=0, y index=1]{data/Pressure-N=0,1.txt};
\addplot[mark = square] table[x index =0, y index=2]{data/Pressure-N=0,1.txt};
\addplot[mark = triangle] table[x index =0, y index=3]{data/Pressure-N=0,1.txt};
\end{axis}
\end{tikzpicture}
\end{center}
\caption{Pressure profile for different values of $M\in\{0,0.5,1\}$, with $N = 0.1$, $R_c = 0.01$, $\nu_b = 0.1$, $\delta = 0.01$.}
\label{Fig:p1}
\end{figure}
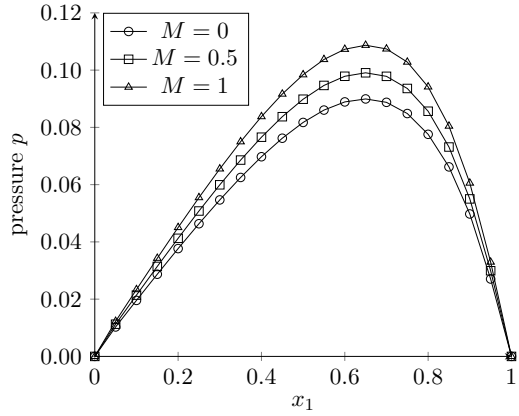

\begin{figure}
\begin{center}
\begin{tikzpicture}[scale=0.8]
\begin{axis}[ axis x line=bottom, axis y line = left,
xlabel={$x_1$}, ylabel={Pressure $p$},
ymin=0, ymax=0.12,
y tick label style={
        /pgf/number format/.cd,
            fixed,
            fixed zerofill,
            precision=2,
        /tikz/.cd
    }
,
legend entries={$M=0$, $M=0.5$, $M=1$}, legend style={at={(0.02,1)}, anchor=north west}]
\addplot[mark = o] table[x index=0, y index=1]{data/Pressure-N=0,2.txt};
\addplot[mark = square] table[x index =0, y index=2]{data/Pressure-N=0,2.txt};
\addplot[mark = triangle] table[x index =0, y index=3]{data/Pressure-N=0,2.txt};
\end{axis}
\end{tikzpicture}
\end{center}
\caption{Pressure profile for different values of $M\in\{0,0.5,1\}$, with $N = 0.2$, $R_c = 0.01$, $\nu_b = 0.1$, $\delta = 0.01$.}
\label{Fig:p2}
\end{figure}
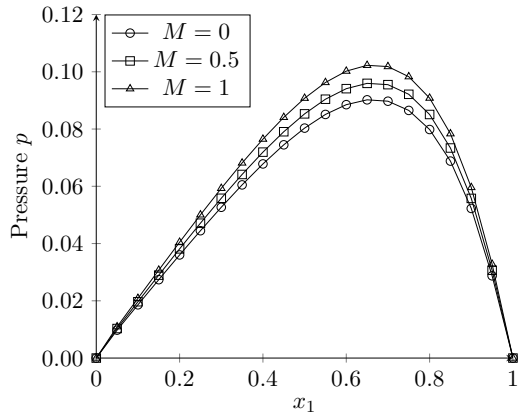

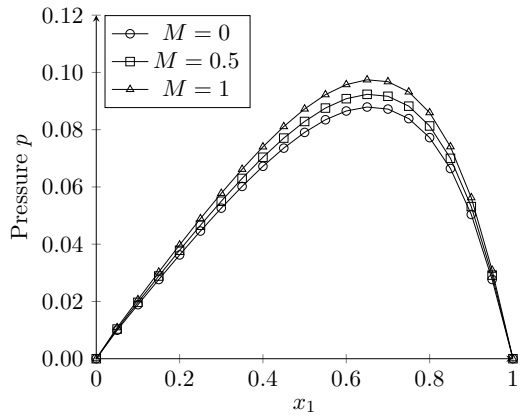
\begin{figure}
\begin{center}
\begin{tikzpicture}[scale=0.8]
\begin{axis}[ axis x line=bottom, axis y line = left,
xlabel={$x_1$}, ylabel={Pressure $p$},
ymin=0, ymax=0.12,
y tick label style={
        /pgf/number format/.cd,
            fixed,
            fixed zerofill,
            precision=2,
        /tikz/.cd
    }
,
legend entries={$M=0$, $M=0.5$, $M=1$}, legend style={at={(0.02,1)}, anchor=north west}]
\addplot[mark = o] table[x index=0, y index=1]{data/Pressure-N=0,3.txt};
\addplot[mark = square] table[x index =0, y index=2]{data/Pressure-N=0,3.txt};
\addplot[mark = triangle] table[x index =0, y index=3]{data/Pressure-N=0,3.txt};
\end{axis}
\end{tikzpicture}
\end{center}
\caption{Pressure profile for different values of $M\in\{0,0.5,1\}$, with $N = 0.3$, $R_c = 0.01$, $\nu_b = 0.1$, $\delta = 0.01$.}
\label{Fig:p3}
\end{figure}


\begin{figure}
\begin{center}
\begin{tikzpicture}[scale=0.8]
\begin{axis}[ axis x line=bottom, axis y line = left,
xlabel={Value of $M$}, ylabel={$W/W_0$},
ymin=0.95, ymax=1.35,
y tick label style={
        /pgf/number format/.cd,
            fixed,
            fixed zerofill,
            precision=2,
        /tikz/.cd
    }
,
legend entries={$N=0.1$, $N=0.2$, $N=0.3$}, legend style={at={(0.02,1)}, anchor=north west}]
\addplot[mark = *] table[x index=0, y index=1]{data/Results-N=0,1.txt};
\addplot[mark = o] table[x index =0, y index=1]{data/Results-N=0,2.txt};
\addplot[mark = +] table[x index =0, y index=1]{data/Results-N=0,3.txt};
\end{axis}
\end{tikzpicture}
\end{center}
\caption{Relative load carrying capacity $W/W_0$ plotted against $M$, for different values of $N\in\{0.1,0.2,0.3\}$ and $R_c = 0.01$, $\nu_b = 0.1$, $\delta = 0.01$.}
\label{Fig:Load}
\end{figure}
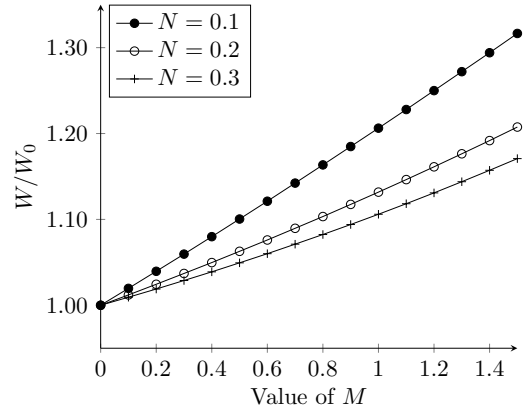




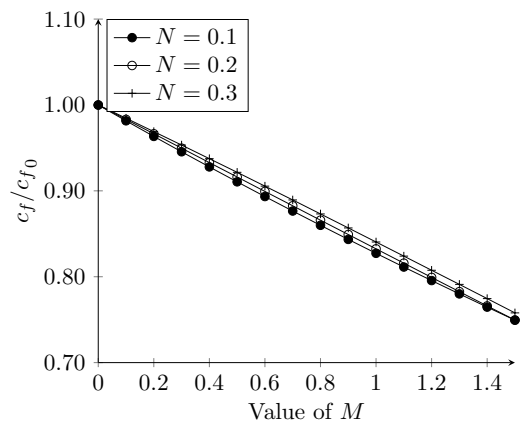
\begin{figure}
\begin{center}
\begin{tikzpicture}[scale=0.8]
\begin{axis}[ axis x line=bottom, axis y line = left,
xlabel={Value of $M$}, ylabel={$c_f/{c_f}_0$},
ymin=0.7, ymax=1.1,
y tick label style={
        /pgf/number format/.cd,
            fixed,
            fixed zerofill,
            precision=2,
        /tikz/.cd
    }
,
legend entries={$N=0.1$, $N=0.2$, $N=0.3$}, legend style={at={(0.02,1)}, anchor=north west}]
\addplot[mark = *] table[x index=0, y index=3]{data/Results-N=0,1.txt};
\addplot[mark = o] table[x index =0, y index=3]{data/Results-N=0,2.txt};
\addplot[mark = +] table[x index =0, y index=3]{data/Results-N=0,3.txt};
\end{axis}
\end{tikzpicture}
\end{center}
\caption{Relative coefficient of friction $c_f/{c_f}_0$ plotted against $M$, for different values of $N\in\{0.1,0.2,0.3\}$ and $R_c = 0.01$, $\nu_b = 0.1$, $\delta = 0.01$.}
\label{Fig:FrictionCoef}
\end{figure}

\newpage
\section*{Acknowledgements}
The second author has been partially supported by the Croatian Science Foundation (Grant No.~3955) and University of Zagreb (Grant No.~202778). The third author has been supported by the projects MTM2011-24457 of the \emph{Ministerio de Econom\'ia y Competitividad} and FQM309 of the \emph{Junta de Andaluc\'ia}.





\section*{References}
\begin{thebibliography}{00}
\bibitem{Szeri} A.Z.~Szeri, Fluid Film Lubrication: Theory and Design, Cambridge University Press, 1998.
\bibitem{Reynolds} O.~Reynolds, On the theory of lubrication and its applications to Mr. Beauchamp Towers'
experiments, including an experimental determination of the
viscosity of olive oil, Philos.~Trans.~Roy.~Soc.~London 177 (1886)
157--234.
\bibitem{Bayada1} G.~Bayada, M.~Chambat, The transition between the Stokes equations and the Reynolds equation: a mathematical proof, Appl.~Math.~Opt.~14 (1986)
73--93.
\bibitem{Archiv} A.~Duvnjak, E.~Maru\v{s}i\'{c}-Paloka, Derivation of the
Reynolds equation for lubrication of a rotating shaft, Archivum
Mat.~(Brno) 36 (2000) 239--253.
\bibitem{Wil} J.~Wilkening, Practical error estimates for Reynolds'
lubrication approximation and its higher order corrections, SIAM J.~
Math.~Anal.~41 (2009) 588--630.
\bibitem{Ja1} E.~Maru\v{s}i\'{c}-Paloka, I.~Pa\v{z}anin, S.~Maru\v{s}i\'{c}, Second order model in fluid film lubrication, C.R.~Mecanique 340 (2012) 596--601.
\bibitem{exp1} G.J.~Jonhnston, R.~Wayte, H.A.~Spikes,
The measurement and study of very thin lubricant films in concentrated contacts, Tribol.~Trans.~\textbf{34} (1991) 187--194.
\bibitem{exp2} J.B.~Luo, P.~Huang, S.Z.~Wen,
Thin film lubrication part I: study on the transition between EHL and thin film lubrication using relative optical interference intensity technique, Wear \textbf{194} (1996) 107--115.
\bibitem{exp3} J.B.~Luo, P.~Huang, S.Z.~Wen, L.~Lawrence,
Characteristics of fluid lubricant films at nano-scale, J.~Tribol.~\textbf{121} (1999) 872--878.
\bibitem{Eringen} A.C.~Eringen, Theory of micropolar fluids, J.~Math.~Mech.~16 (1966) 1--16.
\bibitem{Panasenko1} D.~Dupuy, G.~Panasenko and R.~Stavre,
Asymptotic methods for micropolar fluids in a tube structure,
Math.~Mod.~Meth.~Appl.~Sci.~\textbf{14} (2004) 735--758.
\bibitem{Panasenko2} D.~Dupuy, G.~Panasenko and R.~Stavre, Asymptotic solution for a micropolar flow in a curvilinear channel, Z.~Angew.~Math.~Mech.~\textbf{88}, (2008) 793--807.
\bibitem{Ja2} I.~Pa\v{z}anin, Asymptotic behavior of micropolar fluid flow through a
curved pipe, Acta Appl.~Math.~\textbf{116} (2011) 1--25.
\bibitem{Ja3} I.~Pa\v{z}anin, Modeling of solute dispersion in a circular pipe filled with micropolar fluid, Math.~Comput.~Model.~57 (2012) 2366--2373.
\bibitem{Ja4} M.~Bene\v{s} I.~Pa\v{z}anin, F.~J.~Su\'arez-Grau, Heat flow through a thin cooled pipe filled with a micropolar fluid, J.~Theor.~Appl.~Mech.~53 (2015) 569--579.
\bibitem{Bresch} D.~Bresch, C.~Choquet, L.~Chupin, T.~Colin, M.~Gisclon, Roughness-induced effect at main order on the Reynolds approximation, SIAM Multiscale Model.~Simul.~8 (2010) 997--1017.
\bibitem{fiz} J-L.~Ligier. Lubrification des paliers moteurs. Technip, 1997.
\bibitem{Bayada2} G.~Bayada, M.~Chambat, New models in the theory of the hydrodynamic lubrication of rough surfaces, J.~Tribol.~110 (1988)
402--407.
\bibitem{Benh} N.~Benhaboucha, M.~Chambat, I.~Ciuperca, Asymptotic behaviour of pressure and stresses in a thin film flow with a rough boundary, Quart.~Appl.~Math.~63 (2005)
369--400.
\bibitem{Chupin} L.~Chupin, S.~Martin, Rigorous derivation of the thin film approximation with roughness-induced correctors, SIAM J.~Math.~Anal.~44 (2012) 3041--3070.
\bibitem{Ja5} I.~Pa\v{z}anin, F.~J.~Su\'arez-Grau, Effects of rough boundary on the heat transfer in a thin-film flow, C.~R.~Mecanique 341 (2013) 646--652.
\bibitem{Sinha2} J.~Prakash, P.~Sinha, Lubrication theory for micropolar fluids and its application to a journal bearing, Int.~J.~Engng.~Sci.~13 (1975) 217--232.
\bibitem{Bayada3} G.~Bayada, G.~Lukaszewicz, On micropolar fluids in the theory of lubrication. Rigorous derivation of an analogue of the Reynolds equation,
Int.~J.~Engng.~Sci.~34 (1996) 1477--1490.
\bibitem{Ja6} E.~Maru\v{s}i\'{c}-Paloka, I.~Pa\v{z}anin, S.~Maru\v{s}i\'{c}, An effective model for the lubrication with micropolar fluid,
Mech.~Res.~Commun.~\textbf{52} (2013) 69--73.
\bibitem{Boukrouche} M.~Boukrouche, L.~Paoli, Asymptotic analysis of a micropolar fluid flow in a thin domain with a free and rough boundary, SIAM J.~Math.~Anal.~44 (2012) 1211-1256.
\bibitem{Ja7} I.~Pa\v{z}anin, F.~J.~Su\'arez-Grau, Analysis of the thin film flow in a rough thin domain filled with micropolar fluid, Comput.~Math.~Appl.~68 (2014) 1915--1932.
\bibitem{Bessonov1} N.M. Bessonov, Boundary viscosity conception in hydrodynamical theory of lubrication. Russian Academy of Sicence, Institute of the Problems of Mechanical Engineering, St-Petersbourg, preprint 81 (1993) 105--108.
\bibitem{Bessonov2} N.M. Bessonov, A new generalization of the Reynolds equation for a micropolar fluid and its application to bearing theory, Tribol. Int. 27 (1994) 105-108.
\bibitem{Bayada4} G.~Bayada, N.~Benhaboucha, M.~Chambat, New models in micropolar fluid and their application to lubrication, Math. Mod. Meth. Appl. Sci. 15  (2005) 343--374.
\bibitem{Ja8} I.~Pa\v{z}anin, Asymptotic analysis of the lubrication problem with nonstandard boundary conditions for microrotation, to appear in FILOMAT (2016).
\bibitem{Lukas} G.~Lukaszewicz, Micropolar Fluids: Theory and Applications, Birkh\"{a}user, Boston, 1999.
\bibitem{Ja9} I.~Pa\v{z}anin, On the micropolar flow in a circular pipe: the effects of the viscosity coefficients, Theor.~Appl.~Mech.~Lett.~1 (2011) 062004.
\bibitem{Allaire} G.~Allaire, Homogenization and two-scale convergence, SIAM J.~Math.~Anal.~23 (1992) 1482--1518.

\end{thebibliography}
\end{document}